\documentclass[twocolumn]{autart} 

\usepackage{amsmath,amssymb,amsfonts}

\usepackage{dsfont}
\usepackage{mathrsfs}
\usepackage{mathtools}

\usepackage{algorithm, algorithmicx, algpseudocode}

\makeatletter
\renewcommand{\@biblabel}[1]{}
\makeatother

\begin{document}

\begin{frontmatter}

\title{ A Dual Koopman Approach to Observer
Design \\ for Nonlinear Systems}

  \thanks[footnoteinfo]{Part of these results were presented at the 26th International Symposium on Mathematical Theory of Networks and Systems by J.~Mohet. Corresponding author J. Mohet (Tel. +32 81 724 929.).}

  \author[Namur]{Judicael Mohet\thanksref{footnoteinfo}}\ead{judicael.mohet@unamur.be}
 \author[Namur]{Alexandre Mauroy}\ead{alexandre.mauroy@unamur.be}
  \author[Namur]{Joseph J.J. Winkin}\ead{joseph.winkin@unamur.be}

\address[Namur]{University of Namur, Department of Mathematics and Namur Research Institute for Complex Systems (naXys), Rue de Bruxelles 61, B-5000 Namur, Belgium}   
          
\begin{keyword}                            
    nonlinear observer and filter design; Koopman operator; analytic design; infinite-dimensional systems; Hardy space on the polydisc; observability and detectability.
\end{keyword}

\begin{abstract}   
The Koopman operator approach to the state estimation problem for nonlinear systems is a promising research area. The main goal of this paper is an attempt to provide a rigorous theoretical framework for this approach. In particular, the (linear) dual Koopman system is introduced and studied in an infinite dimensional context. Moreover, new concepts of observability and detectability are defined in the dual Koopman system, which are shown to be equivalent to the observability and detectability of the nonlinear system, respectively. The theoretical framework is applied to a class of holomorphic dynamics. For this class, a Luenberger-type observer is designed for the dual Koopman system via a spectral method, yielding an estimate of the state of the nonlinear system. A particular attention is given to the existence of an appropriate solution to the dual Koopman system and  observer, which are defined in the Hardy space on the polydisc. Spectral observability and detectability conditions are derived in this setting, and the exponential convergence of the Koopman observer is shown. Finally, numerical experiments support the theoretical findings.
\end{abstract}

\end{frontmatter}

\section{Introduction}

State estimation is a prerequisite to monitor a system and therefore appears as one of the main topics in control theory. In the case of linear systems, the popular Luenberger observer and the Kalman filter allow asymptotic state estimation with an exponentially decaying error. However, state estimation for nonlinear systems is much more involved, mostly due to the fact that controlling the nonlinear error dynamics is a nontrivial task. A partial solution to the nonlinear state estimation problem relies on a local linearization of the dynamics. For instance, the celebrated Extended Kalman filter is based on this approach, where a Kalman filter is applied on the nonlinear dynamics, by solving the Riccati equation based on the linearized dynamics and the output map. The main drawback is that only local convergence is guaranteed, apart from specific structured dynamics.

Besides techniques based on linearization, a variety of methods have been developed for nonlinear estimation (see \cite{ber22} for an overview). Some methods are based on the differential detectability condition (\cite[Definition 3.2]{ber22}). However, this condition implies the existence of an observer that converges only if the initial error is small enough. Better convergence results can be obtained by imposing additional conditions which are quite restrictive. We refer the reader to e.g. \cite{san11}, \cite{san15}, \cite{pra01}, and \cite{man18} for more details. Another approach is to use a Kalman-like observer for state-affine systems in normal form, see e.g. \cite{ham90}. In this case, turning the original system into a state-affine system in normal form is key to observe it. A few works in this direction are presented in \cite{kre83}, \cite{bos89}, and \cite{jou03}, where necessary and sufficient conditions are given for the existence of a linearizing change of variables making the resulting linear system observable. In \cite{And06}, a sufficient condition for the existence of a generalized Luenberger observer is given. Moreover, it is proved that a linearizing change of coordinate is the solution of a partial differential equation. Nevertheless, finding an appropriate change of variables remains a challenge, thereby limiting the applicability of the methods. 

In regard to the above-mentioned challenges of nonlinear state estimation, the Koopman operator theory appears to be an appealing approach. Indeed, the Koopman operator allows describing nonlinear dynamics through linear infinite-dimensional systems, and approximated finite-dimensional ones, for which classical linear estimation techniques can be directly used. 
To the authors' knowledge, the Koopman operator framework was applied to the state estimation problem for the first time in \cite{sur16-2} and \cite{sur16}. In the first work, the authors define a Koopman observer form based on spectral properties of the Koopman operator. This form allows estimating the state of a discrete-time nonlinear system linearly. In \cite{sur16}, the Koopman observer form is generalized to input-output continuous-time systems and rank conditions are given in order to characterize the observability of the nonlinear dynamics.

Koopman operator-based approach to observability was also studied in \cite{yeu17} through the Gramian operator and applied to biological systems in \cite{bal24}. Moreover, a sufficient condition for the (approximated) observability of the infinite dimensional Koopman system was given in \cite{mes21}, based on the symmetries of the nonlinear dynamics. Finally, an observer based on the Koopman operator was designed in \cite{ram21} for nonlinear systems over finite fields in the context of cryptography, and the recent works \cite{jad23}, \cite{jun23} and \cite{dah24} further developed the design of Koopman operator based observers in various data-driven situations.

The recent literature clearly illustrates that the observability problem and the design of observers based on the Koopman operator framework is an attractive and promising approach. However, this body of work relies on finite-dimensional spectral representations, assuming that the output map lies in a finite-dimensional span of the eigenfunctions of the Koopman operator, or equivalently, that an invariant linear subspace is known, which contains the output map. Such an assumption is quite difficult to check, and is even not satisfied in general, but relies on the idea that the approximation error vanishes as the dimension of the subspace grows to infinity. In this context, the main goal of this paper is to support the development of past and future Koopman operator based observers, by constructing a theoretical framework in infinite dimension, without any approximation. To do so, we reformulate the state estimation problem in the dual Koopman system framework, where we pinpoint specific properties that are equivalent to observability and detectability of the nonlinear system. In the specific case of a stable hyperbolic equilibrium, we design an infinite-dimensional Luenberger observer (see e.g. \cite{gre75}) for the dual Koopman system defined in the Hardy space on the polydisc. We show that the estimation error converges arbitrarily exponentially fast to zero, and we validate this theoretical framework with numerical experiments.

This paper is organized as follows. In Section \ref{section:Koopman}, basic definitions and properties in Koopman operator theory are recalled and the dual Koopman system and observer are defined. Next, the observability and detectability analysis is performed in Section \ref{section:obs} and, in particular, the equivalence between the observability of the nonlinear system and the observability of the dual Koopman system is proved. The case of the Hardy space on the polydisc is analyzed in Section \ref{section:hardy}, where the well-posedness of the dual Koopman system is shown. In Section \ref{section:KO}, convergence results of an infinite-dimensional Luenberger observer and observability/detectability criteria are obtained. Numerical examples are presented in Section \ref{section:Num} and concluding remarks are given in Section \ref{Conclusion}.

\section{Koopman operator theory and the dual Koopman system} \label{section:Koopman}
A reminder of the Koopman operator theory is provided in the first part of this section. Next, the dual Koopman system and observer are introduced. The section ends with the definition of a concept of solution adapted to these systems.

\subsection{Preliminaries on Koopman operator theory}
Let us consider a nonlinear dynamical system on $\mathbb R^n$ governed by the differential equation
\begin{equation} \label{nonlinODE_no_output}
            \dot{\mathbf x}(t) = F(\mathbf x(t)),  \quad\mathbf x(0)=\mathbf{x}_0.
 \end{equation}
Under regularity assumptions on the vector field $F$ (e.g. $F\in C^1(\mathbb R^n)$), the system admits a unique solution given by the flow $\varphi^t(\mathbf{x}_0)$, where
$\varphi^{\cdot}\colon \mathbb R^+ \to C^1(\mathbb R^n)$ is a semigroup, i.e. $\varphi^0 = I$, where $I$ is the identity and $\varphi^{t+s} = \varphi^t \circ \varphi^s,$ for all $t,s > 0$.
 The Koopman operator associated with the system \eqref{nonlinODE} is the composition operator
\begin{equation*}
    K(t) \colon X \to X, f \mapsto K(t)f =  f \circ \varphi^t, \forall t \geqslant 0,
\end{equation*}
where $X$ is a Banach space.
Moreover, the Koopman system on $X$ associated with \eqref{nonlinODE_no_output} is defined as
\begin{equation}\label{koopmansystem}
  \dot f(t) = A_F f(t),\quad f(0) = f_0,      
\end{equation}
where $f_0, f(t)\in X$ are the observables and $A_F\colon D(A_F) \subset X \to X$ is the infinitesimal generator of the Koopman semigroup. When it exists, the solution is given by $f(t) = K(t) f_0 = f_0 \circ \varphi^t$. 
This implies that there is a one-to-one correspondence between the Koopman system \eqref{koopmansystem} and the dynamical system \eqref{nonlinODE_no_output}. In particular, the knowledge of the Koopman operator provides complete information on the evolution of the state of \eqref{nonlinODE_no_output} through the observables  $f_0$. For example, picking $f_0$ as the $i$th component of the identity function, i.e. $f_0(\mathbf{x})=x_i$, leads to
\begin{equation}\label{koopamsol}
    f(t)(\mathbf{x}_0) = K(t) f_0(\mathbf{x}_0) = (\varphi^t)_i(\mathbf{x}_0), ~\forall t \geqslant 0.
\end{equation}
If the dynamics $F$ is smooth enough, the operator $A_F$ is the (unbounded) linear operator 
\begin{equation*}\label{infinitesimal_gene}
        A_Ff = F \cdot \nabla f = \sum_{k=1}^n F_k \frac{\partial f}{\partial x_k} 
\end{equation*}
defined on the domain 
\begin{equation*}
    D(A_F) = \left\{ f\in X \colon F \cdot \nabla f \in X \right\}.
\end{equation*}
The main advantage of the above approach is that the Koopman system \eqref{koopmansystem} is linear, even though the dynamics $F$ is not.
It follows that spectral properties of the Koopman system are well-defined and can be leveraged in well-chosen spaces $X$ for various purposes, such as the design of a linear observer for the nonlinear dynamics \eqref{nonlinODE_no_output}. We refer the interested reader to the textbook \cite{mau20} on Koopman operator in system and control theory.

In this work, we will consider a \emph{reproducing kernel Hilbert space} (RKHS) $X$, i.e. a Hilbert space where every evaluation functional is bounded and therefore represented by a \textit{reproducing function} $k_{\mathbf{x}}$ such that
\begin{equation}
\label{eq:repr_prop}
    \langle f, k_{\mathbf{x}} \rangle = f(\mathbf{x}), \forall f\in X,
\end{equation}
where $\langle \cdot, \cdot \rangle$ is the inner product on $X$. In this context, it is well-known that 
\begin{equation}
\label{eq:repr_prop2}
    K(t)^* k_{\mathbf{x}} = k_{\varphi^t(\mathbf{x})},
\end{equation}
where $K(t)^*$ is the adjoint operator of $K(t)$. In this paper, we will focus on the Hardy space on the polydisc, which is a RKHS (see Sections \ref{section:hardy}, \ref{section:KO} and \ref{section:Num}). This will provide spectral properties of the operator that are useful in our context of state estimation.

\subsection{Dual Koopman system and dual Koopman observer}
Our main goal is to design a state space observer for the nonlinear system
\begin{equation} \label{nonlinODE}
    \left\{  
        \begin{array}{l}
            \dot{\mathbf x}(t) = F(\mathbf x(t)), \quad \mathbf x(0)=\mathbf{x}_0,\\
            \mathbf y(t) = h(\mathbf x(t)),
        \end{array}
    \right.
 \end{equation}
 where $h\colon \mathbb R^n \to \mathbb R^m$ is a possibly nonlinear output map. Most observers are of the form
\begin{equation}
\label{obs_finite}
   \left\{ 
       \begin{array}{l}
           \dot{\hat {\mathbf x}}(t) = F(\mathbf{\hat x}(t)) + L(\mathbf{\hat y}(t)-\mathbf y(t)), \quad \mathbf{\hat x}(0)=\mathbf{\hat{x}}_0\\
           \mathbf{\hat y}(t) = h(\mathbf{\hat x}(t)),
       \end{array}
   \right.
\end{equation}
where $\mathbf{\hat x}$ is the estimate of the state.
When $F$, $h$ and $L$ are linear, \eqref{obs_finite} is the well-known Luenberger observer. In this case, the error dynamics is purely linear and  can be controlled through the output error injection linear map $L$. However, if $F$ or $h$ is nonlinear, the error dynamics is also nonlinear, an issue which can be overcome by linearizing $F$ and $h$ and by designing the Luenberger observer for the approximated linear system. But obviously, the estimation technique is valid only locally in such a case.

\begin{figure}
    \includegraphics[width = \linewidth]{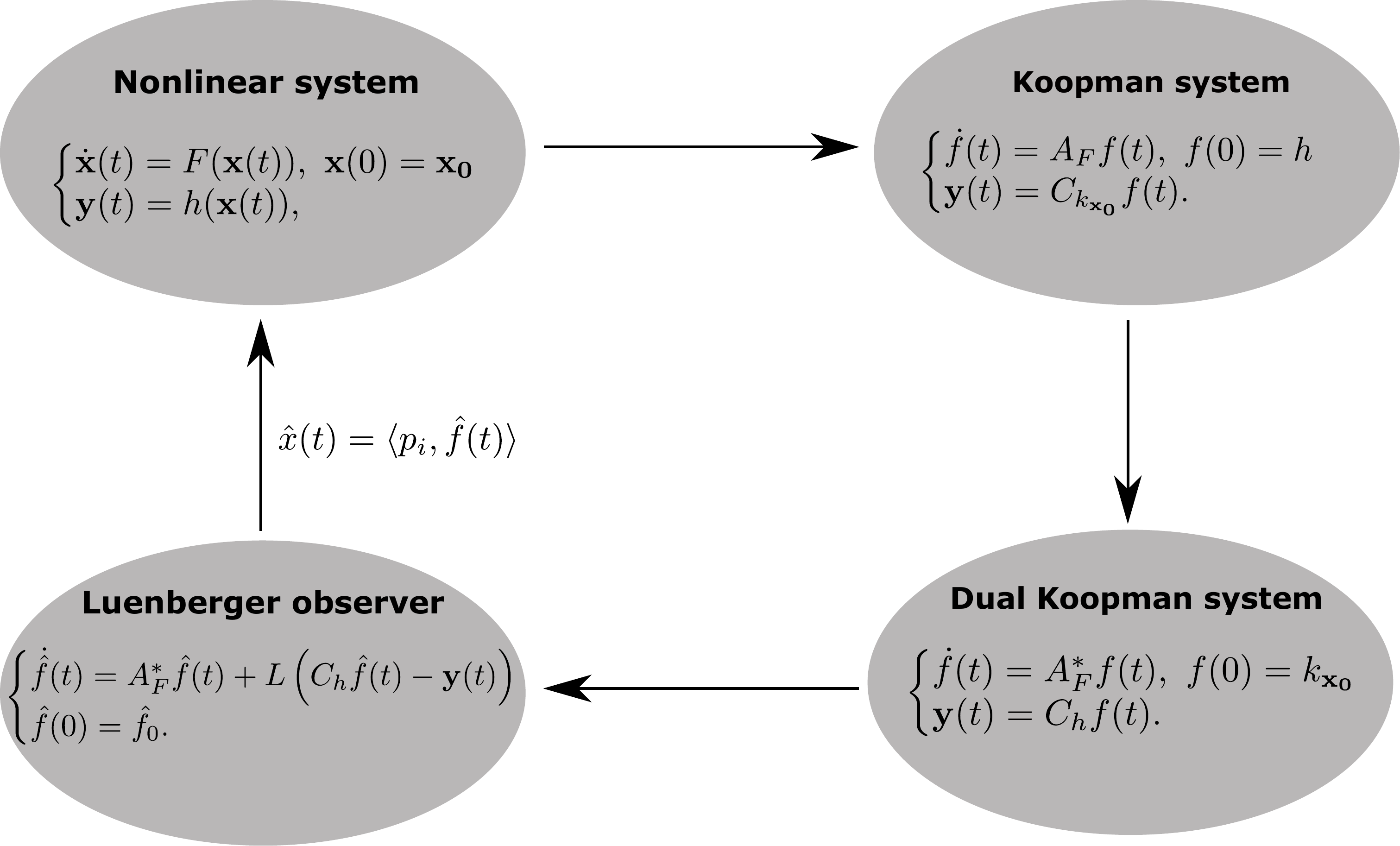}  
    \caption{State estimation with the dual Koopman observer.}
    \label{schema}
\end{figure}

This is where the Koopman operator comes into play. Instead of directly estimating the state $\mathbf{x}(t)$ of the nonlinear dynamical system \eqref{nonlinODE}, we could aim at estimating the observable $f(t)$ of the equivalent linear Koopman system 
\begin{equation}\label{primalkoopmansys}
    \left\{ \begin{array}{l}
        \dot f(t) = A_F f(t), \quad f(0)= h \\
       \mathbf y(t) = C_{\mathbf{x}_0} f(t),
    \end{array}
    \right.
\end{equation}
with the output operator $C_{\mathbf{x}_0} f := f(\mathbf{x}_0)$, so that we verify that $C_{\mathbf{x}_0} f(t) = (K(t) h)(\mathbf{x}_0) = h(\varphi^t(\mathbf{x}_0))=y(t)$. However, the output operator depends on the initial condition and is therefore unknown, while the initial state $f_0=h$ is known. Note also that, if $f_0$ is chosen to be the identity function instead of $h$, then it follows from \eqref{koopamsol} that estimating $f(t)$ is equivalent to estimating the values of the flow at \emph{all} initial conditions, but not at the actual initial state, which remains unknown. For these reasons, the Koopman system \eqref{primalkoopmansys} does not seem appropriate in the context of state estimation. Instead, we will rely on the following dual formulation.
\begin{defn}
     Let $X$ be a reproducing kernel Hilbert space defined on $S\subseteq \mathbb C^n$. The dual Koopman system of \eqref{nonlinODE} is given by
\begin{equation}\label{dualkoopmansys}
    \left\{ \begin{array}{l}
        \dot f(t) = A_F^* f(t), \quad f(0)= k_{\mathbf{x}_0} \\
       \mathbf y(t) = C_h f(t),
    \end{array}
    \right.
\end{equation}
where $f(t) \in X$ for all $t\geqslant 0$, $A_F^*$ is the adjoint operator of $A_F$ and the output operator $C_h$ is defined by
\begin{equation*}
    C_h g = \begin{pmatrix}
        \langle g,{h_1}\rangle & \ldots 
& \langle g, h_m \rangle
    \end{pmatrix}^T, ~ g\in X,
\end{equation*}
where $h = (h_1, \cdots,h_m) \in X^m$.
 \end{defn}
Note that the RKHS framework is needed to properly define the reproducing function $k_{\mathbf{x}_0}$. 
Moreover, by using the reproducing properties \eqref{eq:repr_prop}-\eqref{eq:repr_prop2}, one gets 
$$ C_h f(t)= C_h K(t)^* k_{\mathbf{x}_0}  = C_h k_{\varphi^t(\mathbf{x}_0)} = h(\varphi^t(\mathbf{x}_0)) = \mathbf y(t).$$ 
The dual formulation is appropriate to state estimation. Actually, the output map is known, while the state $f(t)=k_{\varphi^t(\mathbf{x}_0)}$ is directly related to the unknown trajectory $\varphi^t(\mathbf{x}_0)$. This implies that state estimation for the nonlinear dynamics \eqref{nonlinODE} is equivalent to state estimation for the dual Koopman system \eqref{dualkoopmansys}. Indeed, if $\hat f(t)$ is an estimate of $k_{\varphi^t(\mathbf{x}_0)}$ and if $p_k(\mathbf z) = z_k$ is the $k$th component of the identity function, then $\langle p_k, \hat f(t)\rangle$ is an estimate of 
\begin{equation} \label{nonlinearstate}
    \langle p_k,k_{\varphi^t(\mathbf{x}_0)}\rangle = (\varphi^t(\mathbf{x}_0))_k.
\end{equation}
Moreover, since system (\ref{dualkoopmansys}) is linear, a Luenberger observer seems natural. It is provided by a (bounded linear) output injection operator $L$ such that the spectrum of $A^*_F + L C_{h}$ lies in the left-half complex plane.
The existence of such an operator is investigated in the next section. In this case, the Koopman observer (KO) is given by 
\begin{equation} \label{Koopman_obs}
    \dot{\hat f}(t) = A_F^* \hat f(t) + L \left( C_h \hat f(t) - \mathbf y(t)\right), \quad \hat f(0)=\hat f_0
\end{equation}
and the dynamics of the error $\varepsilon = \hat f - f$ is governed by
\begin{equation*}
    \dot{\varepsilon}(t) = (A_F^*+LC_h) \varepsilon(t), \quad \varepsilon(0) = \hat f_0 - f_0.
\end{equation*}
Our methodology is summarized in Figure \ref{schema}.

\subsection{Pseudo-weak solutions}

A classical solution to the dual Koopman system \eqref{dualkoopmansys}, starting from any initial condition, can be obtained if the adjoint Koopman semigroup is strongly continuous. However, we do not need this property here: the dual Koopman system is only defined with a reproducing function as initial condition. Moreover, we only have to estimate \eqref{nonlinearstate} so that we can focus on the solution in a weak sense. These two observations naturally lead to the concept of pseudo-weak solution.

\begin{defn}\label{def:weaksol}
Let $X$ be a RKHS defined on $S\subseteq \mathbb C^n$ and $A$ be a closed and densely defined linear operator on $X$. Consider the abstract Cauchy problem 
\begin{equation} \label{eqgeneral}
    \dot f(t) = A f(t) + \Gamma(t),\quad f(0) =f_0,
\end{equation}
where $\Gamma(t) \in X$ for all $\; t\geqslant 0$ and $f_0 \in \text{span}\left\{ k_{\mathbf x} \colon \mathbf x \in S\right\}$.
A \emph{pseudo-weak solution} to \eqref{eqgeneral} is a function $t\mapsto f(t)$ such that
 \begin{enumerate}
     \item $f\in C([0,\tau],X)$ for all $\tau >0$, i.e. $f(t)$ is strongly continuous on $[0,\tau],$
     \item $\langle f(t), g\rangle$ is absolutely continuous for all $g\in D(A^*)$,
     \item for all $g\in D(A^*)$, 
     \begin{equation*}
        \frac{d}{dt} \langle f(t),g\rangle = \langle f(t), A^*g\rangle + \langle \Gamma(t),g\rangle.
     \end{equation*}
 \end{enumerate}

    \end{defn}
Under appropriate regularity conditions on the function $\Gamma(\cdot)$, if $A$ is the generator of a strongly continuous semigroup on $X$, then a pseudo-weak solution coincides with a mild or classical solution. 
The existence and uniqueness of a pseudo-weak solution for the dual Koopman system and the dual Koopman observer defined on the Hardy space on the polydisc will be proven in Proposition \ref{solution} and Theorem \ref{obs_convergence}, respectively.

\section{Observability and detectability}\label{section:obs}

In addition to the fact that the dual Koopman system \eqref{dualkoopmansys} is directly related to the finite dimensional nonlinear system \eqref{nonlinODE}, its observability and detectability properties are also directly connected to those of the finite dimensional system. This connection is investigated in this section. In what follows, we suppose that $K(t)$ and $K(t)^*$ are bounded on $X$ and that \eqref{dualkoopmansys} and \eqref{Koopman_obs} admit unique pseudo-weak solutions on $X$.

We first recall the concept of observability of the finite dimensional nonlinear system \eqref{nonlinODE} on a subset $S\subseteq \mathbb C^n$.
\begin{defn}
The (finite dimensional) nonlinear system \eqref{nonlinODE} is said to be
\begin{enumerate}
   \item observable in finite time $\tau>0$ if, for all $\mathbf z, \mathbf{z'} \in S,$
   $$h(\varphi^t(\mathbf z)) =  h(\varphi^t(\mathbf z')), \forall t \in [0,\tau]$$
   implies that $\mathbf z = \mathbf z'.$
   \item observable in infinite time if, for all $\mathbf z, \mathbf{z'} \in S,$
   $$h(\varphi^t(\mathbf z)) =  h(\varphi^t(\mathbf z')), \forall t \geqslant 0$$
   implies that $\mathbf z = \mathbf z'$
\end{enumerate} 
\end{defn}

Similar concepts are defined for infinite dimensional systems, based on the classical notion of 
approximate observability (see e.g. \cite[Definition 6.2.12]{zwa20}). In the specific context of the dual Koopman system, we propose the weaker notion of pointwise approximate observability.
\begin{defn}
    Consider the dual Koopman system \eqref{dualkoopmansys} on the RKHS $X$ defined on $S\subseteq \mathbb C^n$. For any positive $\tau$, define the observability map
    \begin{equation*}
       \mathcal C^\tau\colon X \to L^2([0,\tau],\mathbb R^m), f \mapsto C_h K^*(\cdot) f,
    \end{equation*}
    and consider the set 
    $\mathcal K = \{ k_{\mathbf z} - k_{\mathbf{z'}}\colon \mathbf z, \mathbf{z'} \in S\}$.
    The dual Koopman system is said to be
    \begin{enumerate}
        \item pointwise approximately observable (PAO) in finite time $\tau > 0$ if the restriction of the observability map to the set $\mathcal K$ is injective, i.e.
        \begin{equation*}
            \ker \mathcal C^\tau\vert_{\mathcal{K}} = \{ 0\}.
        \end{equation*}
        \item pointwise approximately observable in infinite time if 
        \begin{equation*}
            \bigcap_{\tau >0} \ker \mathcal C^\tau\vert_{\mathcal K} = \{ 0\}.
        \end{equation*}
    \end{enumerate}
\end{defn}

We are now in a position to establish the equivalence between observability of the finite dimensional nonlinear system and pointwise approximate observability of the dual Koopman system. We note that approximate observability is a stronger property implying that the observability map is injective on the whole space $X$, therefore it is \textit{not} equivalent to observability of the associated finite-dimensional dynamics.
\begin{prop}\label{linkPAOobser}
    Consider the dual Koopman system \eqref{dualkoopmansys}, associated with the nonlinear system \eqref{nonlinODE} on the RKHS $X$.
    Then, the dual Koopman system is pointwise approximately observable in infinite time (finite time $\tau >0$, respectively) 
    if and only if the finite dimensional system \eqref{nonlinODE} is observable in infinite time (finite time $\tau >0$, respectively).
\end{prop}

\begin{pf}
The result directly follows from the equality
\begin{equation*}
    \mathcal C^\tau\left( k_{\mathbf{z}} - k_{\mathbf{z'}}\right) = h(\varphi^{\cdot}(\mathbf z)) - h(\varphi^\cdot(\mathbf{z'}))
\end{equation*}
and the fact that the reproducing functions $k_\mathbf{z}$ and $k_{\mathbf{z'}}$ are equal if and only if $\mathbf z=\mathbf{z'}$.
\hfill\ensuremath{\blacksquare} \end{pf}

When the nonlinear system \eqref{nonlinODE} is not observable, one can rely on asymptotic detectability.
\begin{defn}
The finite dimensional system \eqref{nonlinODE} is $\beta-$asymptotically detectable, with $\beta <0$, if, for all $\gamma <-\beta$ and for all $\mathbf z, \mathbf{z'} \in S$, 
$$h(\varphi^t(\mathbf z)) = h(\varphi^t(\mathbf{z'})) \quad \forall t>0 $$
implies that 
$$\lim_{t\to \infty}e^{\gamma t}\left\vert \varphi_k^t(\mathbf z) - \varphi_k^t(\mathbf{z'})\right\vert =0, \;  \forall k=1,\ldots,n. $$ 
\end{defn}
Let us now define the corresponding concept for the dual Koopman system.
\begin{defn}
The dual Koopman system \eqref{dualkoopmansys} is pointwise $\beta-$detectable, with $\beta <0$, if for all $\gamma <- \beta$ and for all $\mathbf z, \mathbf{z'} \in S$,
$$k_\mathbf{z} -k_{\mathbf{z'}} \in \displaystyle \bigcap_{\tau > 0}\ker \mathcal C^\tau $$
implies that
$$\lim_{t\to\infty}e^{\gamma t}\vert \langle p_k, K(t)^*(k_\mathbf{z}-k_{\mathbf{z'}})\rangle \vert = 0, \quad \forall k=1,\ldots,n,
$$
where $p_k(\mathbf z) = z_k$.
\end{defn}

Again, the following result shows the equivalence between the above definitions related to detectability.
\begin{prop}
    The dual Koopman system \eqref{dualkoopmansys} is pointwise $\beta-$detectable, for some $\beta<0$, if and only if the associated 
    finite dimensional system \eqref{nonlinODE} is $\beta-$asymptotically detectable.
\end{prop}
\begin{pf}
    Observe that, for all $k=1,\ldots,n$,
    \begin{align*}
        \vert \langle p_k, K(t)^*(k_\mathbf{z}-k_{\mathbf{z'}})\rangle \vert& = \vert \langle K(t)p_k, k_\mathbf{z}-k_{\mathbf{z'}}\rangle\vert\\
        &=\vert \varphi^t_k(\mathbf z) - \varphi^t_k(\mathbf{z'}) \vert.
    \end{align*}
\hfill\ensuremath{\blacksquare} \end{pf}
For finite dimensional nonlinear systems, asymptotic detectability is a necessary condition for the existence of an 
(asymptotic) observer \cite[Section 3]{ber22}. This is also true for the dual Koopman observer.
\begin{prop}\label{link_cvg_det}
Suppose that the Koopman observer converges weakly and exponentially fast at rate $\beta <0$, that is,
for all $\mathbf z, \hat{\mathbf{z}}\in S $ and for all $g\in D(A_F)$, 
\begin{equation*} 
    e^{\gamma t} \left\vert \langle g, \hat f(t) -  f(t)\rangle  \right\vert \xrightarrow{t\to\infty} 0 \quad \forall \gamma < -\beta,
\end{equation*}
where $\hat f(t)$ is the solution to the Koopman observer \eqref{Koopman_obs} with $\hat f_0 = k_{\hat{\mathbf{z}}}$ and $f(t)$ is the solution 
to the dual Koopman system \eqref{dualkoopmansys} with $f_0 = k_{\mathbf z}$. Then the dual Koopman system is pointwise $\beta-$detectable.
\end{prop}
\begin{pf}
Pick any $\mathbf z$ and $\mathbf{z'} \in S$ such that $k_{\mathbf z} -k_{\mathbf{z'}}\in \displaystyle\bigcap_{\tau >0 }\mathcal \ker \, C^\tau$. Then
 \begin{equation*}
    \mathbf{y}(t) = C_h K(t)^* k_{\mathbf z} = C_h K(t)^* k_{\mathbf{z'}}.
 \end{equation*}
 It follows that the solution $\hat f(t)$ to the Koopman observer \eqref{Koopman_obs} is identical for the initial states $k_{\mathbf z}$ and $k_{\mathbf{z'}}$.
Hence, with $f(t) := K(t)^*k_z$ and $f'(t) := K(t)^*k_{z'}$, for all $k=1,\ldots,n$ and for all $\gamma < -\beta$, 
\begin{align*}
    e^{\gamma t} \left\vert \langle p_k, f(t) -f'(t)\rangle  \right\vert &\leqslant e^{\gamma t}\left\vert \langle p_k, f(t) -\hat f(t)\rangle \right\vert\\
     &+ e^{\gamma t}  \left\vert \langle p_k, \hat f(t) -f'(t)\rangle  \right\vert,
    \end{align*}
    where the right-hand side tends to $0$ as $t$ tends to infinity by assumption.
\hfill\ensuremath{\blacksquare} \end{pf}

\section{The dual Koopman system on the Hardy space on the polydisc} \label{section:hardy}

From this point on, we will focus on the Hardy space on the unit polydisc. When defined in this space of observables and associated with a specific class of dynamics, the Koopman operator possesses spectral properties that are useful to state observer design. Also, the existence and uniqueness of pseudo-weak solutions (Definition \ref{def:weaksol}) for the dual Koopman system and observer hold in this context.

\subsection{The Hardy space on the polydisc}

The unit polydisc ${\mathbb D}^n$ is the Cartesian product of $n$ unit discs of $\mathbb C$, i.e.
\begin{equation*}
    {\mathbb D}^n = \{ \mathbf{z}\in \mathbb C^n \colon \forall k=1,\ldots, n, \vert z_k\vert < 1\}.
\end{equation*}
Its topological boundary is $\partial\left(\mathbb D^n\right) = \{\mathbf z \in \mathbb C^n \colon \exists k\in\{1,\ldots,n\}, \vert z_k \vert = 1\}$, and its distinguishable boundary is 
$\mathbb {\mathbb T}^n = \left(\partial \mathbb D\right)^n =  \{ \mathbf{z}\in \mathbb C^n \colon \forall k=1,\ldots,n, \vert z_k\vert =1 \}$.
The Hardy space on ${\mathbb D}^n$, denoted as $\mathbb H^2({\mathbb D}^n)$, is the space of all holomorphic functions $f\colon {\mathbb D}^n\to \mathbb C$ such that the norm 
\begin{equation*}
    \Vert f \Vert^2 =\lim_{r\to 1^-} \int_{\mathbb {\mathbb T}^n}\vert f\left(r\mathbf{w}\right)\vert^2dm_n(\mathbf{w})
\end{equation*}
is finite, where $m_n$ is the normalized Lebesgue measure\footnote{i.e. if $\mathbf w = (e^{i\theta_1}, \ldots,e^{i\theta_n})$, then $\int_{\mathbb {\mathbb T}^n} f(r\mathbf{w}) dm_n(\mathbf{w})=$
\begin{equation*}
\frac{1}{(2\pi)^n} \int_0^{2\pi}\cdots \int_0^{2\pi} f(r e^{i\theta_1},\ldots,re^{i\theta_n}) d\theta_1\cdots d\theta_n.
\end{equation*}} on $\mathbb {\mathbb T}^n$. It is a Hilbert space equipped with the inner product 
\begin{equation*}\label{innerH2}
\left\langle f,g\right\rangle=\lim_{r\to 1^-}\int_{ \mathbb T^n}f\left(r\mathbf{w}\right) \bar g\left(r\mathbf{w}\right) dm_n(\mathbf{w}).
\end{equation*} 
Moreover, this space is a reproducing kernel Hilbert space with the reproducing functions
\begin{equation*}\label{kernelH2}
k_{\mathbf{z_0}}(z) = \displaystyle \prod_{i=1}^n \dfrac{1}{1-\overline{{z_0}_i} z_i},
\end{equation*}
see for example \cite{pau16}.

In what follows, $\mathbb N$ denotes the set of integers including zero.
The Hardy space $\mathbb H^2(\mathbb{D}^n)$ admits the orthonormal basis of monomials $\{e_\alpha\}_{\alpha \in \mathbb N^n}$, with $e_\alpha(\mathbf z) = \mathbf{z}^\alpha=z_1^{\alpha_1} ~ \cdots~ z_n^{\alpha_n}$.
The total degree of $e_\alpha$ will be denoted by \mbox{$\vert \alpha \vert = \alpha_1 + \cdots + \alpha_n$}.
We will assume without loss of generality that the monomials are ordered according to the lexicographic order.\footnote{The lexicographic order relation is defined by 
    $\alpha  < \beta \Leftrightarrow \vert \alpha \vert < \vert \beta \vert$  or 
    $\vert \alpha \vert = \vert \beta \vert \wedge (\alpha_j > \beta_j, j=\min\{i \colon \alpha_i \neq \beta_i\})$
for all $\alpha, \beta \in \mathbb N^n$.} Note that the monomials $e_\alpha$ with $|\alpha|=1$ are the projections $p_k$. In addition, $\mathbb H^2({\mathbb D}^n)$ is a separable Hilbert space so that there is an isomorphism between $\mathbb H^2({\mathbb D}^n)$ and $l^2(\mathbb N^n)$, the space of square-summable sequences in $\mathbb C$. Specifically, since $f\in \mathbb H^2({\mathbb D}^n)$ is analytic on ${\mathbb D}^n$, there is a unique sequence $f_\alpha \in l^2(\mathbb N^n)$, 
defined by the Taylor coefficients $f_\alpha = \langle f, e_\alpha \rangle$ of $f$ such that $f = \sum_{\alpha\in \mathbb N^n} f_\alpha e_\alpha$. 
As a result, the norm and the inner product can be rewritten as
\begin{equation*}
   \Vert f\Vert^2 = \Vert (f_\alpha)_{\alpha\in \mathbb N^n} \Vert_{2}^2
\end{equation*}
and
\begin{equation}
\label{eq:inner_prod}
    \langle f,g\rangle= \langle (f_\alpha)_{\alpha\in \mathbb N^n} , (g_\alpha)_{\alpha\in \mathbb N^n} \rangle_2 = \sum_{\alpha \in \mathbb N^n} f_\alpha \overline{g_{\alpha}}.
\end{equation}
In this context, the operators $A_F$ and $A_F^*$ are both closed and densely defined, see e.g. \cite{ros24}. However, the boundedness of the Koopman operator and its adjoint is not guaranteed. We will prove this property under additional assumptions in the next section.

\subsection{Well-posedness of the dual Koopman system on $\mathbb H^2(\mathbb{D}^n)$}

In this section, we provide a class of dynamics such that the boundedness of the Koopman operator and its adjoint holds on $\mathbb H^2(\mathbb{D}^n)$. The existence of a pseudo-weak solution to the dual Koopman system \eqref{dualkoopmansys} is proved in this case. The class of dynamics satisfies the following assumptions.
\begin{itemize}
    \item[] \textbf{Assumption 1 (holomorphic flow).} The vector field $F\colon \mathbb D^n\cup \mathbb T^n \to \mathbb C^n$ is holomorphic on $\mathbb D^n$ and continuous on $\mathbb T^n$, and generates a flow 
    $$ \mathbb R^+ \times \mathbb D^n \to \mathbb D^n, (t,\mathbf z) \mapsto \varphi^t(\mathbf z)$$
     which, for all $t\geqslant 0$, is holomorphic on $\mathbb D^n$.\\
    
    \item[] \textbf{Assumption 2 (forward invariance).} For all $t> 0$, there exists $R_t \in (0,1)$ such that
    \begin{equation*}
     \forall \mathbf z \in \mathbb T^n, \varphi^t(\mathbf z) \in D(\mathbf 0,R_t)^n.
    \end{equation*}
   where $D(\mathbf a,R)$ denotes the open disc centered at $\mathbf a$ of radius $R$.
\end{itemize}

The following result, based on \cite[Theorem 1.7]{che16}, provide sufficient conditions on the vector field $F$ so that the above assumptions are satisfied.
\begin{prop}\label{chelthm1.7}
    Suppose that $F\colon \overline{\mathbb D^n} \to \mathbb C^n$ is of the form 
    \begin{equation*}\label{holo_flow_cond}
        F_j(\mathbf z) = \left(z_j - u_j(\hat{\mathbf{z}}_j)\right) G_j(\mathbf z), \quad j=1,\ldots,n
    \end{equation*}
    with $\mathbf{\hat z}_j = \begin{pmatrix}
        z_1, \ldots, z_{j-1},z_{j+1}, \ldots, z_n
    \end{pmatrix}$, and where $u_j \colon \overline{\mathbb D^{n-1}} \to \overline{\mathbb D}$ and $G_j \colon \overline{\mathbb D^n}\to \mathbb C$ are holomorphic on $\mathbb D^n$ and satisfy
    \begin{equation} \label{holo_flow_ineq}
        \text{Re } \left((1-\overline{z_j} u_j(\hat{\mathbf{z}}_j))G_j(\mathbf z) \right) < 0.
    \end{equation}
     Under these conditions, Assumptions 1 and 2 hold.
 \end{prop}
 The proof can be found in the appendix.
 
Thanks to Assumptions 1 and 2, we can prove our first main result on the boundedness of the Koopman operator and its adjoint operator on $\mathbb H^2(\mathbb{D}^n)$. Note that a different proof is also given in \cite{mug24}.
Our proof is based on the following lemma.
\begin{lem}\cite[Theorem 9]{jaf90} \label{jafari1}
    Let $p \in (1,\infty)$. The composition operator $f\mapsto f\circ \phi$ is bounded on $\mathbb H^p({\mathbb D}^n)$ if and only if 
    \begin{equation*}
        \label{jafari_bounded_cond}
        \sup_{\mathbf{z}\in {\mathbb D}^n} \displaystyle \int_{{\mathbb T}^n} \prod_{j=1}^n \dfrac{1-\vert z_j\vert^2}{\vert 1 - \overline{z_j}\phi^*_j\vert^2} dm_n < +\infty,
    \end{equation*}
where $\displaystyle \phi^*(\mathbf w) = \lim_{r\to1^-} \phi(r\mathbf w), \mathbf w\in {\mathbb T}^n$. 
\end{lem}

\begin{prop}
Consider the dynamical system \eqref{nonlinODE} that satisfies Assumptions 1 and 2. Then, the adjoint of the Koopman operator $K(t)^*$ associated with the dynamics is bounded on $\mathbb H^2(\mathbb{D}^n)$.
  \end{prop}
\begin{pf}
Note that Assumption 1 implies that the composition operator associated with $\varphi^t$ is well-defined on $\mathbb H^2({\mathbb D}^n)$. Moreover, the continuity of the flow on $\mathbb T^n$ implies that 
$(\varphi^t)^* = \varphi^t$, where $(\varphi^t)^*(\mathbf{w}) = \lim_{r\to 1^-} \varphi^t(r\mathbf{w}), \mathbf{w} \in {\mathbb T}^n$.
By Assumption 2, for all $t > 0$, there exists $R_t \in (0,1)$ such that
\begin{equation*}
    \vert 1 - \overline{z}_j\varphi_j^t(\mathbf w) \vert \geqslant \vert 1 - \vert \overline{z}_j\varphi^t_j(\mathbf w) \vert~\vert
    >  1 - R_t
\end{equation*}
for all $\mathbf z\in {\mathbb D}^n$, $\mathbf w\in {\mathbb T}^n$ and $j=1,\ldots,n$.
Then,
\begin{align*}
    &\displaystyle \int_{{\mathbb T}^n} \prod_{j=1}^n \dfrac{1 - \vert z_j\vert^2}{\vert 1 - \overline{z}_j (\varphi^t(\mathbf w))\vert^2} dm_n(\mathbf w) \\
    < ~& \frac{1}{\left( 1 - R_t\right)^{2n}} \prod_{j=1}^n (1 - \vert z_j\vert^2)\\
    <~& \frac{1}{\left( 1 - R_t\right)^{2n}}
\end{align*}
for all $\mathbf z\in {\mathbb D}^n$. The boundedness of the Koopman operator follows from the last inequality, by Lemma \ref{jafari1}. Hence, its adjoint is also bounded and $\Vert K(t)^*\Vert = \Vert K(t) \Vert$ for all $t \geqslant 0$. 
\hfill\ensuremath{\blacksquare} \end{pf}

In our context, the adjoint Koopman semigroup is strongly continuous (and invariant) in the dense subspace spanned by the reproducing functions $k_\mathbf{z}, \mathbf z \in \mathbb D^n$. Hence, by \cite[Proposition 5.3]{eng01}, strong continuity holds if and only if the operator-valued function $K(\cdot)^*$ is bounded on some interval $[0,\delta]$. This property is still an open question. However, we can prove the existence of a pseudo-weak solution (see Definition \ref{def:weaksol}).

\begin{prop}  \label{solution}
    Suppose that Assumptions 1 and 2 hold. Then, the dual Koopman system \eqref{dualkoopmansys} admits a unique pseudo-weak solution on $\mathbb H^2(\mathbb{D}^n)$, which is given by 
    \begin{equation*}
           f(t) = K(t)^*k_{\mathbf{x}_0} = k_{\varphi^t(\mathbf{x}_0)},~t\geqslant 0.
    \end{equation*}
\end{prop}

\begin{pf}
First, the map $f\colon [0,+\infty) \to \mathbb H^2(\mathbb{D}^n)$ is continuous, since $(K(t)^*)_{t\geqslant 0}$ is strongly continuous on the set of reproducing functions. Secondly, the function $t\mapsto\varphi^t(\mathbf{x}_0)$ is continuously differentiable. Hence, for all $g\in D(A_F)$, the map $t\mapsto \langle f(t),g\rangle = g(\varphi^t(\mathbf{x}_0))$ is absolutely continuous. Finally, observe that 
\begin{equation*}
    \frac{d}{dt} \langle f(t),g\rangle = \nabla g \cdot F(\varphi^t(\mathbf{x}_0))= \langle f(t), A_fg\rangle
\end{equation*}
for all $g\in D(A_F)$.
This implies that $f(t)$ is a pseudo-weak solution to \eqref{dualkoopmansys}. The proof of uniqueness is standard and similar to the proof given in \cite{ball77}. 
\hfill\ensuremath{\blacksquare}\end{pf}

\subsection{Spectral analysis}

We can now derive specific spectral properties of the Koopman operator that we will leverage for state estimation. Those properties require the following additional assumptions on the nonlinear system \eqref{nonlinODE}.
\begin{itemize}
    \item[] \textbf{Assumption 3 (stable hyperbolic equilibrium).} The dynamics admits a stable equilibrium at the origin, whose basin of attraction contains the polydisc $\overline{{\mathbb D}^n}$. Moreover, the eigenvalues $\lambda_j$ of the Jacobian matrix $\mathcal J_F(\mathbf 0)$ are simple and satisfy $\text{Re}(\mathbf \lambda_j) < 0$ for all $j=1,\dots,n$.\\
    \item[] \textbf{Assumption 4 (non-resonant eigenvalues).} The eigenvalues $\lambda_j$ of the Jacobian matrix $\mathcal J_F(\mathbf 0)$ are nonresonant, i.e. 
    for all $(m_1,\ldots,m_n) \in \mathbb{N}^n$ satisfying $ \sum_{l=1}^n m_l\geqslant 2$, 
    $$ \lambda_j \neq \sum_{l=1}^n m_l \lambda_l, ~\forall j=1,\ldots,n.$$
    \item[] \textbf{Assumption 5 (output map).} The components of the output map $h$ belong to $\mathbb H^2(\mathbb{D}^n)$.
\end{itemize}
\begin{rem}

Nonresonant eigenvalues are required to rely on the Poincaré-Dulac linearization theorem. However, other linearization theorems exist, with different assumptions, see e.g. the Siegle-Bruno theorem in \cite{ber23}. In \cite{kre01}, similar assumptions were considered to obtain a necessary and sufficient condition for the existence of a  change of variable that linearizes the dynamics up to a nonlinear injection term.
\end{rem}
The first remarkable  property resulting from those assumptions is the fact that the operators $A_F$ and $A_F^*$ admit a series expansion that allows to represent both operators as infinite matrices.
\begin{lem}
For all $f\in D(A_F)$, 
\begin{equation*}
    A_F f = \sum_{\alpha\in \mathbb N^n} \left( \sum_{\beta\in\mathbb N^n} {A_F}_{\alpha,\beta} f_\beta \right)e_\alpha\\
\end{equation*}
where $f_\beta = \langle f, e_\beta\rangle$ and ${A_F}_{\alpha,\beta} =\langle A_F e_\beta, e_\alpha\rangle$.
\end{lem}
\begin{pf}
Since $\{e_\alpha\}_{\alpha\in\mathbb N^n}$ is an orthonormal basis of $\mathbb H^2(\mathbb{D}^n)$, any $f\in D(A_F)$ can be expanded as $f = \sum_{\beta\in \mathbb N^n} f_{\beta} e_\beta$. Hence, we can prove that
\begin{equation} \label{toprove}
 A_F \left( \sum_{\beta\in \mathbb N^n} f_{\beta} e_\beta \right) = \sum_{\beta\in \mathbb N^n} f_{\beta} A_Fe_\beta.
\end{equation}
Indeed, in \cite{mug23}, it is shown that the right-hand side is given by
\begin{equation}\label{mugAebeta}
\begin{split}
   & \sum_{\beta\in \mathbb N^n} f_{\beta} A_Fe_\beta \\
   & \quad = \sum_{l=1}^n F_l \sum_{\beta \in \mathbb N^n} (\beta_l+1) f_{(\beta_1, \ldots, \beta_{l-1}, \beta_l + 1, \beta_{l+1},\ldots, \beta_n)} e_\beta.
   \end{split}
\end{equation}
It remains to show that the left-hand side of \eqref{toprove} is also equal to the right-hand side of \eqref{mugAebeta}. 
To do so, consider the operators $A_1$ and 
$A_2$ given by 
\begin{equation*}
    A_1 f= \nabla f \text{ and } A_2 \omega = F \cdot \omega,
\end{equation*}
for all $f \in D(A_F)$ and $\omega \in (\text{Hol}(\mathbb D^n))^n$, respectively. 
Observe that $A_F = A_2 A_1$ on $D(A_F)$. 
Moreover, for all $f\in D(A_F)$ and for all $l=1,\ldots,n$,
\begin{equation*}
    \left(A_1 f\right)_l =  \sum_{\beta \in \mathbb N^n} (\beta_l + 1) f_{(\beta_1, \ldots, \beta_{l-1}, \beta_l + 1, \beta_{l+1},\ldots, \beta_n)} e_\beta.
\end{equation*}
Hence, for all $f\in D(A_F)$,
\begin{align*}
    A_F f& = A_2A_1 f \\
    &= \sum_{l=1}^n F_l \sum_{\beta \in \mathbb N^n} (\beta_l + 1) f_{(\beta_1, \ldots, \beta_{l-1}, \beta_l + 1, \beta_{l+1},\ldots, \beta_n)} e_\beta.
\end{align*}
In view of \eqref{mugAebeta}, it follows that identity \eqref{toprove} holds, which implies that
\begin{align*}
    A_Ff &= \sum_{\alpha\in \mathbb N^n} \langle A_Ff, e_\alpha\rangle e_\alpha\\
    & = \sum_{\alpha\in \mathbb N^n} \left \langle  A_F \left( \sum_{\beta\in \mathbb N^n} f_{\beta} e_\beta \right), e_\alpha \right \rangle e_\alpha \\
    &= \sum_{\alpha\in \mathbb N^n} \sum_{\beta\in\mathbb N^n} \langle A_F e_\beta, e_\alpha\rangle \, f_\beta \, e_\alpha.
\end{align*}
\hfill\ensuremath{\blacksquare}
\end{pf}

It is proved in \cite{mug23} that the infinite matrix representation $\mathcal A_F$ of the operator $A_F$ is lower block triangular of the form
\begin{equation}\label{matrix_rep}
    \mathcal A_Ff= \begin{pmatrix}
        [0] & \cdots & & & \\
        [0]& [\mathcal A_{11}] & [0] &\cdots & \\
        [0] & [\mathcal A_{21}] & [\mathcal A_{22}] & [0] & \cdots \\
        \vdots & \vdots & \vdots & \vdots & \ddots 
    \end{pmatrix}\begin{pmatrix}
        [f_0]\\ [f_1] \\ [f_2] \\ \vdots
    \end{pmatrix},
\end{equation} 
where the block $[\mathcal A_{ij}]$ maps the components $[f_j]$ corresponding to the degree $j$ of $f\in D(A_F)$ (in the monomial's basis) to 
the components corresponding to the degree $i$ of $A_Ff$.
The adjoint operator $A_F^*$ admits the same type of representation, but with the transpose conjugate matrix.

Now, we prove that the so-called \emph{principal eigenfunctions} lie in the Hardy space $\mathbb H^2(\mathbb D^n)$.
\begin{prop}\label{principal_ef}
    Suppose that Assumptions 1-5 hold and consider the principal eigenfunctions $\phi_j$ of the Koopman operator, with $j=1,\dots,n$, which satisfy
    \begin{equation*}
        A_F \phi_j = \lambda_j \phi_j \text{ and } K(t) \phi_j = e^{\lambda_j t} \phi_j \; ,
    \end{equation*}
    where $\lambda_j$ are the eigenvalues of the Jacobian matrix $\mathcal J_F(\mathbf 0)$. Then, $\phi_j \in \mathbb H^2(\mathbb D^n)$ for all $j=1,\dots,n$.
\end{prop}
\begin{pf}
    Assumptions 3 and 4 together with Poincaré linearization theorem (e.g. \cite{bel14}) ensure the existence of a holomorphic conjugacy $\nu$ in a neighborhood $\Omega$ of the origin, which maps the 
    nonlinear system \eqref{nonlinODE} into
    \begin{equation*}
        \left\{ 
            \begin{array}{l}
                \dot{\mathbf y}(t) = \mathcal J_F(\mathbf 0) \mathbf y(t)\\
                \mathbf y(0)=\nu(\mathbf{x}_0),
            \end{array}
        \right.
    \end{equation*}
    where $\mathbf y = \nu(\mathbf x)$. It is known that the principal eigenfunctions of $K(t)$ are given by 
    \begin{equation*}
        \phi_j(\mathbf x) = \mathbf{w}_j^* \cdot \nu(\mathbf x), \quad \mathbf x \in \Omega,
    \end{equation*}
    where $\mathbf{w}_j$ is the left eigenvector of $\mathcal J_F(\mathbf 0)$ associated with $\lambda_j$. Those eigenfunctions can be extended to $\overline{\mathbb D^n}$ by defining, for $\mathbf z\in \overline{\mathbb D^n}$
    \begin{equation*}
        \phi_j(\mathbf z) = e^{-\lambda_j t} \phi_j(\varphi^t(\mathbf z)).
    \end{equation*}    
    Since the origin is stable, there exists $t>0$ large enough such that $\varphi^t(\mathbf z)\in \Omega$. In this case, it is clear that the right hand side is holomorphic on ${\mathbb D^n}$ and continuous on $\mathbb T^n$, which implies that $\phi_j \in \mathbb H^2(\mathbb{D}^n).$\hfill\ensuremath{\blacksquare}
\end{pf}
Thanks to the triangular structure of $A_F^*$, useful spectral properties of the dual Koopman operator can be shown.
\begin{prop}\label{point_spectrum}
    Suppose that Assumptions 1-5 hold.  Let us denote the eigenfunctions of $A_F$ and $A_F^*$ by $\{\phi_\alpha\}_{\alpha\in\mathbb N^n}$ and $\{\psi_\alpha\}_{\alpha\in\mathbb N^n}$, respectively. Then, the following assertions hold.
    \begin{enumerate}
        \item The spectrum of the operator $A_F$ is given by 
        \begin{equation*}\label{spectrum}
            \sigma\left(A_F\right) =\left\{ \sum_{i=1}^n \alpha_i \lambda_i \colon \alpha \in\mathbb N^n\right\}, 
        \end{equation*}
        where $\lambda_j$, $j=1,\dots,n$ are the eigenvalues of the Jacobian matrix $\mathcal J_F(\mathbf 0)$.
        \item The sets of eigenfunctions $\{\phi_\alpha\}_{\alpha\in\mathbb N^n}$ and $\{\psi_\alpha\}_{\alpha\in\mathbb N^n}$ are biorthonormal basis of $\mathbb H^2(\mathbb{D}^n)$.
        \item The operators $A_F$ and $A_F^*$ admit the spectral expansions 
\begin{equation} \begin{split} \label{expansions}
      &A_F f = \sum_{\alpha\in\mathbb N^n} \lambda_\alpha \langle  f, \psi_\alpha\rangle \phi_\alpha\\
      &A_F^*g = \sum_{\alpha\in\mathbb N^n} \overline{\lambda_\alpha}\langle g, \phi_\alpha\rangle\psi_\alpha,
\end{split} \end{equation}
for all $f\in D(A_F)$ and $g\in D(A_F^*),$ where $\{\lambda_\alpha\}_{\alpha\in\mathbb N^n} = \sigma\left(A_F\right)$.
Moreover, the Koopman operator and its adjoint admit the spectral expansions
\begin{align*}
    &K(t) f = \sum_{\alpha\in\mathbb N^n} e^{\lambda_\alpha t} \langle f, \psi_\alpha\rangle \phi_\alpha\\
    &K(t)^* f= \sum_{\alpha\in\mathbb N^n} e^{\overline{\lambda_\alpha} t} \langle f, \phi_\alpha\rangle\psi_\alpha,
\end{align*}
for all $f\in \mathbb H^2(\mathbb{D}^n)$ and $t\geqslant 0$.
     \end{enumerate}
    \end{prop}
   Some of these properties are well-known. In particular, it is known that $\phi_\alpha=\prod_{j=1}^n \phi_j^{\alpha_j}$ satisfies the eigenvalue equation with the corresponding eigenvalue \mbox{$\lambda_\alpha = \sum_{j=1}^n \alpha_j \lambda_j$} (note that $\lambda=0$ is a trivial eigenvalue associated with a constant eigenfunction) and $\phi_\alpha$ belongs to $\mathbb H^2(\mathbb D^n)$ if and only if $\alpha \in \mathbb N^n$. A comprehensive proof of Proposition \ref{point_spectrum} is given in the appendix for the sake of completeness.

\begin{rem}
Proposition \ref{point_spectrum} implies that we can write
$g(\mathbf z) = \sum_{\alpha\in\mathbb N^n} \langle g, \psi_\alpha\rangle \phi_\alpha(\mathbf z)$ for all $g\in\mathbb H^2(\mathbb{D}^n)$ and $\mathbf z\in \mathbb D^n$. By identifying this series with the Taylor expansion of $g$ in the variables $\xi = \Phi(\mathbf z) =\left( \phi_1(\mathbf z), \ldots, \phi_n(\mathbf z)\right)$ (where $\phi_j$ are the principal eigenfunctions), we obtain
\begin{equation*}
    \langle g,\psi_\alpha\rangle = \left.\frac{\partial^\alpha \left(g\circ\Phi^{-1}\right)}{\partial \xi^{\vert \alpha \vert}}\right|_{\mathbf{0}} \quad \forall \alpha\in \mathbb N^n.
\end{equation*}
Note also that the \emph{exact} finite expansion of the eigenfunctions $\psi_\alpha$ in the basis of monomials can be obtained by computing the left eigenvectors of the truncated matrix representation \eqref{matrix_rep}. Similarly, a truncated Taylor expansion of $\phi_\alpha$ can be obtained with the right eigenvectors.
\end{rem}

\section{The dual Koopman observer on $\mathbb H^2(\mathbb{D}^n)$} \label{section:KO}

In this section, we provide novel observability and detectability criteria, construct the dual Koopman observer, and prove its convergence in $\mathbb H^2(\mathbb D^n)$. To do so, we will exploit the spectral properties of the Koopman operator and focus on the estimation of the least stable modes.
For $\beta < 0$, the $\beta-$unstable eigenvalues of the Koopman operator are defined as the eigenvalues $\lambda\neq 0$ of $A_F^*$ that satisfy $\text{Re}(\lambda) > \beta$. The number $N_\beta$ of $\beta-$unstable eigenvalues is finite for all $\beta\in \mathbb R_0^-$. 
For the sake of convenience, we will denote those eigenvalues by $\lambda_{\beta,j}^+$, with $j=1,\dots,N_\beta$, and the associated eigenfunctions of $A_F^*$ by $\psi_{\beta,j}^+$. Moreover, denoting $\mathbb N_0^n$ the set of nonzero elements of $\mathbb N^n$, we can define the sets $I_\beta^+ =\{ \alpha\in \mathbb N_0^n \colon \text{Re} (\lambda_\alpha) > \beta\}$ and $I_\beta^- =\mathbb N_0^n \setminus I_\beta^+$, which are associated with the subspaces \mbox{$H^+_\beta \colon =\text{span } \{\psi_\alpha\colon \alpha \in I_\beta^+\}=\text{span }\{\psi^+_{\beta,j}\}_{j=1}^{N_\beta}$} and \mbox{$H^-_\beta \colon = \overline{\text{span }\{\psi_\alpha\colon \alpha \in I_\beta^-\}}$}. We also define the subspace $H^0 = \text{span}\{\psi_0\},$ where $\psi_0$ is the constant function $1$,  so that $\mathbb H^2(\mathbb D^n) = H^0\bigoplus H^+_\beta \bigoplus H^-_\beta.$
In what follows, we will consider the bounded linear spectral projection operators (see e.g. \cite[Theorem 6.1.19]{buh18})
\begin{equation*}
    \begin{array}{l}
        P^0 \colon \mathbb H^2(\mathbb D^n) \to H^0, \quad g\mapsto \langle g,1 \rangle \,1=g(0)1, \\
        P^+\colon \mathbb H^2(\mathbb D^n) \to H^+_\beta, \quad g \mapsto \sum_{\alpha \in I_\beta^+} \langle g,\phi_\alpha\rangle \psi_\alpha, \\
        P^-=I-P^0-P^+. 
    \end{array}
\end{equation*}
Following the same idea as in \cite{gre75}, we will estimate the $\beta-$unstable part of the error dynamics, taking advantage of the finite dimension of $H^+_\beta$.

\subsection{Observability and detectability criteria}

In Section \ref{section:obs}, the concept of pointwise approximate observability (PAO) was introduced for the dual Koopman system.
We proved that this property is equivalent to the observability of the associated nonlinear finite-dimensional
system. But, while PAO is an infinite-dimensional concept, the Koopman framework is developed for a finite-dimensional system, so that 
it seems natural to seek for a finite-dimensional characterization. This can be obtained through the 
principal eigenfunctions (see Proposition \ref{principal_ef}).

\begin{thm}
    Consider the dual Koopman system \eqref{dualkoopmansys} with Assumptions 1-5. If, for all principal eigenfunctions  
    $\psi_j$ ($j\in\{1,\ldots,n\}$) of the operator $A_F^*$, such that $A_F^*\psi_j = \overline{\lambda_j}\psi_j$, there exists $i\in\{1,\cdots,m\}$ such that 
    \begin{equation}\label{cond_principal_obs}
        \langle h_i, \psi_j \rangle\neq 0,
    \end{equation}
    then the dual Koopman system is PAO in infinite time.
\end{thm}
\begin{pf}
Let us take $\mathbf z,\mathbf{z'}\in {\mathbb D}^n$ such that \mbox{$\mathcal{C}^\tau(k_{\mathbf z}-k_{\mathbf{z'}}) =0$}, for all $\tau > 0$. This means that 
\begin{equation*}
    \sum_{\alpha \in \mathbb N^n} e^{\lambda_\alpha t} \langle  h_i, \psi_\alpha\rangle \langle \phi_\alpha, k_{\mathbf z}-k_{\mathbf{z'}}\rangle = 0,
\end{equation*}
for all $i={1,\cdots,m}$ and for all $t\geqslant 0$. 
Using arguments similar to those of the proof of \cite[Theorem 6.3.4]{zwa20},
it follows from the identity above that 
\begin{equation*}
    \langle h_i, \psi_\alpha\rangle \langle \phi_\alpha, k_{\mathbf z}-k_{\mathbf{z'}}\rangle = 0 \quad \forall i\in\{1,\ldots,m\}, \forall \alpha\in\mathbb N^n.
\end{equation*}
Together with \eqref{cond_principal_obs}, this implies that 
\begin{align*}
    \langle \phi_j, k_{\mathbf z}-k_{\mathbf{z'}}\rangle= 0, \forall j=1,\ldots,n,
\end{align*}
or equivalently,
\begin{equation*}
    \phi_j(\mathbf{z}) = \phi_j(\mathbf{z'}), \forall j=1,\ldots,n,
\end{equation*}
where $\phi_j$, j=1,\ldots, n, are the principal eigenfunctions of $A_F$. This implies that $\mathbf z = \mathbf{z'}$ since $\mathbf{z} \mapsto 
    (\phi_1(\mathbf{z}), \dots, \phi_n(\mathbf{z}))$ is a diffeomorphism and therefore $k_\mathbf{z} = k_\mathbf{z'}$.
\hfill\ensuremath{\blacksquare} \end{pf}
It follows from Proposition \ref{linkPAOobser} that the condition \eqref{cond_principal_obs} implies the observability of the nonlinear system \eqref{nonlinODE}. This condition 
is actually related to the observability criterion given in \cite[Theorem 2.3.4 and Corollary 2.3.5]{isi95} in the geometric control context.

We can also obtain a spectral characterization of $\beta-$asymptotic detectability.
\begin{lem}\label{lemma:det_spectral_criteria}
    Consider the dual Koopman system \eqref{dualkoopmansys} with Assumptions 1-5 and $\beta <0$. If, for all $\beta-$unstable principal eigenfunctions $\psi_j$ (i.e. for all $j$ such that $\psi_j\in H_\beta^+$), there exists $i\in\{1,\ldots,m\}$ satisfying
    \begin{equation} \label{det_cond}
        \langle h_i, \psi_{j}\rangle\neq 0,
    \end{equation} 
    then the dual Koopman system is pointwise $\beta-$detectable.
\end{lem}
\begin{pf}
   Take $\mathbf z, \mathbf{z'}\in {\mathbb D}^n$ such that $\mathcal C^\tau(k_\mathbf{z}-k_\mathbf{z'}) =0$ for all $\tau > 0$. This identity implies that, for all $t\geqslant 0$,
    \begin{equation*}
        \sum_{\alpha\in\mathbb N^n} e^{\lambda_{\alpha}t} \langle h_i, \psi_{\alpha}\rangle \langle \phi_{\alpha}, k_{\mathbf z} - k_{\mathbf{z'}}\rangle = 0, \forall i=1,\ldots,m.
    \end{equation*}
     Together with \eqref{det_cond}, this implies that $\phi_j(\mathbf z)=\phi_j(\mathbf z')$ for all $j$ such that $\textrm{Re}(\lambda_j)>\beta$. Note that, if $\alpha \in I_\beta^+$, then $\textrm{Re}(\lambda_j)>\textrm{Re}(\lambda_\alpha)>\beta$ for all $j$ such that $\alpha_j>0$ and therefore $\phi_\alpha=\Pi_j  \phi_j^{\alpha_j}$ satisfies $\phi_\alpha(\mathbf z)=\phi_\alpha(\mathbf z')$ 
     for all $\phi_\alpha\in H^0\bigoplus H^+_\beta$. It follows that, for $k=1,\ldots,n$,
    \begin{equation*}
        \langle p_k,K(t)^*(k_\mathbf{z}-k_{\mathbf{z'}})\rangle = \sum_{\alpha \in I_\beta^-} e^{\overline{\lambda_\alpha} t} \left \langle p_k,\psi_\alpha\rangle \langle k_{\mathbf z} - k_{\mathbf{z'}}, \phi_\alpha \rangle. \right.
    \end{equation*}
    Now, pick any $t_1 > 0$. Then, for all $t>t_1$, we have
      \begin{equation*}
          \begin{split}
          & \vert \langle p_k,K(t)^*(k_\mathbf{z}-k_{\mathbf{z'}})\rangle  \vert \\
              &\quad \leqslant e^{\beta(t-t_1)} \sum_{\alpha \in I_\beta^-} e^{\text{Re}(\lambda_\alpha)\, t_1} \left\vert \langle p_k,\psi_\alpha\rangle \langle \phi_\alpha, k_{\mathbf z} - k_{\mathbf{z'}}\rangle\right\vert.
          \end{split}
      \end{equation*}
      The above series converges since its general term is the product of terms of a geometric sequence \mbox{$e^{\text{Re}(\lambda_\alpha) t_1}=\mathcal O((e^{-t_1})^{\vert\alpha\vert})$} and a sequence $( \langle p_k,\psi_\alpha\rangle \langle \phi_\alpha, k_{\mathbf z} - k_{\mathbf{z'}}\rangle)_\alpha$ that tends to zero (since $\langle p_k,k_{\mathbf z} - k_{\mathbf{z'}}\rangle = \sum_{\alpha \in \mathbb N^n}\langle p_k,\psi_\alpha\rangle \langle \phi_\alpha, k_{\mathbf z} - k_{\mathbf{z'}}\rangle$). Hence, for all $\gamma < - \beta ,$ \\
      $\displaystyle \lim_{t\to \infty} e^{\gamma t} \vert\langle p_k,K(t)^*(k_\mathbf{z}-k_{\mathbf{z'}})\rangle  \vert =0.$ 
    \hfill\ensuremath{\blacksquare} \end{pf}

\subsection{Design and convergence of the dual Koopman observer}
 
We are now in a position to design the dual Koopman observer and to show its convergence in $\mathbb H^2(\mathbb{D}^n)$ for our class of dynamics.
Using \eqref{expansions} and the bi-orthogonality of the basis, one can show that $H^+_\beta$ and $H^-_\beta$ are $A_F^*$-invariant. It follows that it is natural to define the operators \mbox{$(A^*_F)^0 = A_F^*P^0$}, \mbox{$(A^*_F)^+ = A_F^*P^+$}, \mbox{$(A^*_F)^-=A_F^*P^-, C^0 = C_h P^0$,} $C^+=C_h P^+$, and $C^-=C_h P^-$. 
If the pair $(C^+,(A^*_F)^+)$ is observable, then there exists an operator $L^+:\mathbb{C}^m \to H_\beta^+$ which allows defining the injection operator
  $$L =\begin{pmatrix}
      P^0 L \\P^+L \\ P^-L
  \end{pmatrix} =\begin{pmatrix}
 0\\ L^+\\ 0\end{pmatrix}.$$ 
Then the observer dynamics \eqref{Koopman_obs} takes the form
\begin{equation}\label{observer1}
         \begin{pmatrix}
         \dot{\hat f}^0\\
\dot{\hat f}^+ \\
\dot{\hat{f}}^-     
\end{pmatrix}= \begin{pmatrix}0 & 0 & 0\\
         L^+C^0 & M^+ &   L^+C^- \\
         0 & 0  & (A_F^*)^-
         \end{pmatrix}
         \begin{pmatrix}
         \hat f^0\\
\hat f^+ \\
\hat{f}^-     
\end{pmatrix} - \begin{pmatrix}
            0\\ L^+\\ 0
         \end{pmatrix} y(t)
\end{equation}
where ${\hat{f}}^0(t) =P^0{\hat{f}}(t)$, ${\hat{f}}^+(t) = P^+{\hat{f}}(t)$, ${\hat{f}}^-(t) = P^-{\hat{f}}(t)$, and $M^+ = (A_F^*)^+ + L^+ C^+.$
\begin{rem} \label{rem:constant_mode}
    With the initial condition $\hat{f}^0(0)=k_{\mathbf{\hat{x}_0}}$, it is easy to see that $\hat{f}^0(t) = \hat{f}^0(0) = P^0 k_{\mathbf{\hat{x}_0}} = 1$ and we have $C^0 {\hat f}^0(t) = C_h 1 = h(0)$. It follows that
    \begin{equation*}
    \begin{split}
        \dot{\hat f}^+ & = M^+ {\hat f}^+(t) + L^+ C^- {\hat f}^-(t) - L^+ (h(x(t))-h(0))
        \end{split}
    \end{equation*}
    so that the dynamics of ${\hat f}^+$ can be obtained without considering the component in ${H}^0$ (i.e. the constant mode) by subtracting the constant term $h(0)$ from the output map. This will be exploited in the numerical method (see Algorithm \ref{alg:obs}).
\end{rem}

Next, we can prove our main result. We show that there exists a pseudo-weak solution to \eqref{observer1} and we give a condition such that this solution converges (in the sense of Proposition \ref{link_cvg_det}) to the pseudo-weak solution of the dual Koopman system.
\begin{thm}\label{obs_convergence}
Consider the dynamics \eqref{nonlinODE} with initial condition 
$\mathbf{x}_0$, Assumptions 1 to 5, together with the corresponding dual Koopman system \eqref{dualkoopmansys}, and set $\beta<0$. If, for all $\beta$-unstable eigenfunctions $\psi_{\beta,j}^+$, $j\in\{1,\ldots,N_\beta\}$, of the operator $A_F^*$, there exists $i\in\{1,\ldots,m\}$ such that
\begin{equation} \label{cvg_cond}
    \langle h_i, \psi_{\beta,j}^+\rangle  \neq 0,
\end{equation}
then there exists a bounded linear operator $L$ such that, for all $\mathbf{\hat{x}}_0\in \mathbb D^n,$ the observer dynamics \eqref{observer1} with initial condition $k_{\mathbf{\hat{x}}_0}$ admits a unique pseudo-weak solution $\hat f(t)$. Moreover, this solution converges weakly and exponentially fast at rate $\beta$ to the pseudo-weak solution $f(t)$ of the dual Koopman system, that is
\begin{equation*}
    e^{\gamma t} \left\vert \langle g, \hat f(t) -  f(t)\rangle  \right\vert \xrightarrow{t\to\infty} 0 \quad \forall \gamma < -\beta, \, \forall g \in D(A_F).
\end{equation*}
\end{thm}
\begin{pf}
\textbf{Existence and uniqueness of the solution.} Since $H_\beta^+$ is finite-dimensional, the operator $(A^*_F)^+$ (restricted to $H_\beta^+$) can be represented by the matrix $(\mathcal A_F^*)^+  = \text{diag}(\overline{\lambda_{\beta,1}^+}, \ldots, \overline{\lambda_{\beta,N_\beta}^+})\in \mathbb C^{N_\beta \times N_\beta}$ in the basis $\{\psi_{\beta,j}^+\}_{j=1}^{N_\beta}$. Similarly, the operator $C^+$ (restricted to $H_\beta^+$) is represented by the matrix $\mathcal C^+ \in \mathbb C^{m \times N_\beta}$ with entries $(\mathcal C^+)_{ij} = \langle h_i,\psi_{\beta,j}^+\rangle$. Observe that the condition \eqref{cvg_cond} implies that the pair $(\mathcal C^+,(\mathcal A^*_F)^+)$ is observable. 
Indeed, by the Popov-Belevich-Hautus test, the pair $({\mathcal C}^+,({\mathcal A_F}^*)^+)$ is observable if and only if
     \begin{equation*}
        \text{rank} \begin{pmatrix}
            (\mathcal A^*_F)^+ - \overline{\lambda_{\beta,i}^+} I\\
            \mathcal C^+
        \end{pmatrix} = N_\beta \quad \forall i =1,\ldots,N_\beta.
      \end{equation*}
   The first $N_\beta$ rows form a diagonal matrix whose $ith$ diagonal entry is zero (recall that the eigenvalues of $\mathcal A_F^*$ are simple). Thus, the matrix above is full rank if and only if there is at least one non-zero component in the $i$th column of the matrix $\mathcal C^+$, which is equivalent to the condition \eqref{cvg_cond}. The observability property implies that there exists a matrix $\mathcal{L}^+ \in \mathcal{C}^{N_\beta \times m}$ such that \mbox{$(\mathcal{A}_F^*)^+ + \mathcal{L}^+\mathcal{C}^+$} is $\beta$-stable (e.g. all its eigenvalues have a real part smaller than $\beta$). Note that the matrix $\mathcal{L}^+$ represents $L^+$ in the basis $\{\psi_{\beta,j}^+\}_{j=1}^{N_\beta}$.

Next, denoting $\varepsilon(t)= \hat f(t) - f(t)$, the state estimation error, the abstract Cauchy problem
\begin{equation}\label{errordyn}
    \left\{ \begin{array}{l}
   \dot{\varepsilon}^0(t) = 0,\\
      \dot{\varepsilon}^+(t) = \left( (A^*_F)^+ + L^+ C^+ \right) \varepsilon^+(t) +  L^+C^-  \varepsilon^-(t), \\
         \dot{\varepsilon}^-(t) = (A^*_F)^- \varepsilon^-(t),\\
         \varepsilon(0) = k_{\mathbf{\hat x}_0}-k_{\mathbf x_0},
    \end{array} 
        \right.
\end{equation}
where $\varepsilon^0(t) = P^0\varepsilon(t)  = 0$, and where $\varepsilon^+(t) = P^+\varepsilon(t)$, $\varepsilon^-(t) = P^-\varepsilon(t)$, $\varepsilon(0)^+ = P^+(k_{\mathbf{\hat x}_0}-k_{\mathbf x_0})$, and $\varepsilon(0)^-=P^-(k_{\mathbf{\hat x}_0}-k_{\mathbf x_0})$, admits a pseudo-weak solution on $\mathbb H^2(\mathbb{D}^n)$. Indeed, first observe that both operators $(A_F^*)^+$  and $(A_F^*)^-$ are closed and densely defined. 
Since the dynamics of $\varepsilon^+(t)$ is finite-dimensional, it also admits a pseudo-weak solution on $H^+_\beta$. Thus, it remains to be shown that the $\beta-$stable part admits a pseudo-weak solution on $H^-_\beta$.
It is easy to see that $\varepsilon^-(t)=P^-K(t)^*(k_{\mathbf{\hat x}_0}-k_{\mathbf x_0})$ and to check that the map $t\mapsto \varepsilon^-(t)$ is strongly continuous on $[0,\tau]$ for all $\tau > 0$ since so is $t\mapsto K(t)^*k_{\mathbf{z}}$ and the operator $P^-$ is bounded. Moreover, for all $g\in D((A_F^*)^-)$, 
\begin{align*}
    \langle \varepsilon^-(t),g\rangle &= \langle k_{\varphi^t(\mathbf{\hat{x}}_0)}-k_{\varphi^t(\mathbf{x}_0)},(P^-)^*g\rangle \\
    &=(P^-)^*(g(\varphi^t(\hat{\mathbf x}_0)-g(\varphi^t(\mathbf{x}_0))),
\end{align*}
which is absolutely continuous since $(P^-)^*g \in \mathbb H^2(\mathbb D^n)$. The above relation also implies that 
\begin{align*}
\frac{d}{dt} \langle \varepsilon^-(t),g\rangle & = \nabla(P^-)^*g(\varphi^t(\hat{\mathbf{x}}_0)) \cdot  F(\varphi^t(\hat{\mathbf{x}}_0)) \\
& \quad -\nabla(P^-)^*g(\varphi^t({\mathbf{x}}_0)) \cdot  F(\varphi^t({\mathbf{x}}_0))\\
      &=\langle  k_{\varphi^t(\mathbf{\hat{x}}_0)}-k_{\varphi^t(\mathbf{x}_0)},A_F(P^-)^*g\rangle.
             \end{align*}
It is easy to show that $(P^-)^*g = \displaystyle \sum_{\alpha\in I^-_\beta} \langle g,\psi_\alpha\rangle \phi_\alpha$ so that
     \begin{equation*}
     A_F (P^-)^* g = \displaystyle \sum_{\alpha\in I_\beta^-} \lambda_\alpha \langle g, \psi_\alpha\rangle \phi_\alpha = (P^-)^* A_F  g.
     \end{equation*}
       Hence, we get
     \begin{align*}
        \frac{d}{dt} \langle \varepsilon^-(t),g\rangle &= \langle k_{\varphi^t(\mathbf{\hat{x}}_0)}-k_{\varphi^t(\mathbf{x}_0)}, (P^-)^*A_F g\rangle \\
        &=\langle \varepsilon^-(t),A_Fg\rangle.
     \end{align*}
      This leads to the existence of a pseudo-weak solution for the error dynamics, and for the dual Koopman observer as well.
The uniqueness of this pseudo-weak solution can be established by using arguments of  \cite{ball77}. 

      \textbf{Convergence of the solution.} Pick any $t_1 > 0$. For all $g\in D(A_F)$ and all $t > t_1$, the following inequality holds:
      \begin{equation*}
          \begin{split}
              \vert \langle \varepsilon^-(t),g\rangle  \vert &\leqslant \sum_{\alpha \in I_\beta^-} e^{\text{Re}(\lambda_\alpha)\,(t-t_1)} e^{\text{Re}(\lambda_\alpha)\,t_1} \left\vert \langle \varepsilon(0),\phi_\alpha\rangle \langle \psi_\alpha,g\rangle\right\vert\\
              &\leqslant e^{\beta(t-t_1)} \sum_{\alpha \in I_\beta^-} e^{\text{Re}(\lambda_\alpha)\, t_1} \left\vert \langle \varepsilon(0),\phi_\alpha\rangle \langle \psi_\alpha,g\rangle\right\vert.
          \end{split}
      \end{equation*}
      Using the same arguments as in the proof of Lemma \ref{lemma:det_spectral_criteria}, we can show that $\varepsilon^-(t)$ converges weakly and exponentially fast to $0$.
      Moreover, the solution $\varepsilon^0(t)=0$ does not affect the convergence rate.
      It remains to show the convergence of $\varepsilon^+(t)$. It follows from \eqref{errordyn} that
      \begin{equation*}
          \mathbf{\varepsilon}^+(t) = \Theta(t) \, \mathbf{\varepsilon}^+(0) + (\Theta*\mathcal{L}^+C^-\varepsilon^-)(t),
      \end{equation*}
      where $\mathbf{\varepsilon}^+(t)$ is the vector of components of $e^{+}(t)$ in the basis $\{\psi_{j}\}_{j=1}^{N_\beta}$, the operation $*$ is the convolution product, and $\Theta(t)=e^{((\mathcal{A}_F^*)^+ + \mathcal L^+ \mathcal C^+)t}$. Since \mbox{$(\mathcal{A}_F^*)^+ + \mathcal L^+ \mathcal C^+$} is a $\beta$-stable matrix and $\varepsilon^-(t)$ converges weakly and exponentially fast to zero at rate $\beta$, we obtain that $\mathbf{\varepsilon}^+(t)$ converges exponentially fast to zero at the same rate. \hfill\ensuremath{\blacksquare} \end{pf}

      Although the observer solution is not a reproducing function, even if so is the initial condition $\hat{f}(0)=k_{\mathbf{x}_0}$, Theorem \ref{obs_convergence} implies that it converges weakly and exponentially fast to $f(t)=k_{\varphi^t(\mathbf{x}_0)}$. In addition, if we consider $g$ to be the monomial $e_\alpha=p_k \in D(A_F)$ of degree $1$, with $\alpha_k=1$ and $\alpha_i=0$ for all $i\neq k$, we obtain that $\langle e_\alpha, \hat{f}(t)\rangle - \left(\varphi^t(\mathbf x_0)\right)_k = \langle e_\alpha, \hat{f}(t) - f(t)\rangle  \rightarrow 0$, as $t$ tends to infinity. This provides a linear state space observer for the nonlinear system \eqref{nonlinODE}.
      
\begin{rem}
The condition \eqref{cvg_cond} for the convergence of the Koopman observer is stronger than the sufficient condition \eqref{det_cond} for the $\beta$-pointwise detectability of the dual Koopman system. This is not a contradiction since the convergence of the observer is not equivalent to, but implies $\beta$-pointwise detectability (see Proposition \ref{link_cvg_det}). Note also that, in case the condition \eqref{cvg_cond} is not satisfied, the output can be virtually augmented by considering additional lifted output map components. For instance, if the output map is linear, the condition \eqref{cvg_cond} will not be satisfied (even with full state measurements) for $\beta$-unstable eigenfunctions that have typically no first order Taylor component. In this case, considering additional input map components of the form $(h(x))^\alpha=y^\alpha$ can make condition \eqref{cvg_cond} satisfied.
\end{rem}

\section{Numerical method and examples} \label{section:Num}

This section is devoted to the numerical illustration of our method. We first provide the algorithm describing our Koopman observer based method. Second, we define a particular class of dynamics that satisfies Assumptions 1 to 5. Finally, numerical experiments are given and shown to support our theoretical results.

\subsection{Numerical method}

The numerical method associated with the dual Koopman operator-based observer is summarized in Algorithm \ref{alg:obs}. It relies on a finite-dimensional approximation $\mathcal{A}^*_F$ of  $A_F^*$ in a basis of monomials truncated to some total degree $d$. The expansion of the eigenfunctions $\phi_\alpha$ and $\psi_\alpha$ in the monomial basis (e.g. Taylor coefficients) are obtained as the components of the left and right eigenvectors of $\mathcal{A}^*_F$, respectively (see e.g. \cite{mug23}). The inner product \eqref{eq:inner_prod} is computed with the truncated basis. Note also that the constant zeroth-order mode is neglected, since the value $h(0)$ is subtracted from the output $y(t)$ (see Remark \ref{rem:constant_mode}).

\begin{algorithm}[ht]
    \caption{Dual Koopman Observer}
    \label{alg:obs}
    \begin{algorithmic}[1]
    \Statex{\bf Input:} Initial guess $\mathbf{\hat{x}}_0$, truncation degree $d$, Taylor coefficients $F_{l,\gamma}$ with $F_l=\sum_\gamma F_{l,\gamma} \, e_\gamma$ ($l=1,\dots,n$), Taylor coefficients $h_{j,\gamma}$ with $h_j=\sum_\gamma h_{j,\gamma} \, e_\gamma$ \mbox{($j=1,\dots,m$).}
    \Statex{\bf Output:} Estimated state $\mathbf{\hat{x}}(t)$.
\State Compute the matrix approximation $[\mathcal{A}^*_F]_d\in\mathbb C^{N_d\times N_d}$ (with $N_d = \frac{(n+d)!}{n!~d!}-1$) of $A_F^*$ with entries
      $  ({\mathcal{A}^*_F})_{\alpha,\gamma} = \sum_{l=1}^n \gamma_l F_{l,(\alpha -\gamma)_l}$
    if $d \geq\vert \alpha \vert\geqslant \vert \gamma \vert>0$, where $$(\alpha -\gamma)_l = (\alpha_1 - \gamma_1, \ldots,\alpha_l - \gamma_l +1, \ldots, \alpha_n -\gamma_n).$$
\State Compute the eigenvalues $\lambda_\alpha$, left-eigenvectors $\mathbf v_\alpha$, and right-eigenvectors $\mathbf w_\alpha$ of the matrix $[\mathcal{A}^*_F]_d$, with $(\mathbf w_\alpha)^*\mathbf v_\alpha=1$ for all $\alpha$.
\State Compute the approximations $\langle h_j,\psi_\alpha \rangle \approx (\mathbf w_\alpha)^* \mathbf{h}_j$, where $\mathbf{h}_j$ is the vector containing the Taylor coefficients $h_{j,\gamma}$, with $0<|\gamma| \leq d$.
\State Choose $\beta < 0$ such that \eqref{cvg_cond} holds, 
and define a bijective map $\zeta: \mathbb{N} \to \mathbb{N}^n$ such that $\zeta(k) \in I_\beta^+$ for all $k \in \{1,\dots,N_\beta\}$.
\State Construct the matrices $(\mathcal A_F^*)^+  = \text{diag}\{\overline{\lambda_{\zeta(k)}}\}_{k=1}^{N_\beta}$, $(\mathcal A_F^*)^-  = \text{diag}\{\overline{\lambda_{\zeta(k)}}\}_{k=N_\beta+1}^{N_d}$, and $\mathcal C^+$, $\mathcal C^-$ with entries $(\mathbf v_{\zeta(j)})^* \mathbf{h}_i$ (with $j=1,\dots,N_\beta$ for $\mathcal C^+$ and \mbox{$j=N_\beta+1,\dots, N_d$} for $\mathcal C^-$).
\State Compute the matrix $\mathcal L^+$ such that all the eigenvalues of \mbox{$(\mathcal{A}_F^*)^+ + \mathcal L^+ \mathcal C^+$} are $\beta$-stable (using a standard pole placement method).
    \State Compute the solution $\hat{\mathbf f}(t)$ to 
    \begin{equation*}
        \dot{\hat{\mathbf f}} = \begin{pmatrix}
         (\mathcal A^*_F)^+ + \mathcal L^+ \mathcal C^+) &   \mathcal L^+\mathcal C^-\\
         0 & (\mathcal A_F^*)^-         
         \end{pmatrix} \hat{\mathbf f}
        - \begin{pmatrix}
             \mathcal L^+\\ 0
         \end{pmatrix} (y(t)-h(0))
    \end{equation*}
        with the initial condition $\mathbf{\hat{f}}_0 = \mathcal V^* \mathbf{e}(\mathbf{\hat{x}}_0),$ where $\mathcal V$ is the matrix with columns $\mathbf v_\alpha$ and $\mathbf{e}(\mathbf{\hat{x}}_0)$ is the vector of components $e_\gamma(\mathbf{\hat{x}}_0)$, $|\gamma|\leq d$.
    \State For $k=1,\dots,n$, recover the state estimate
    \begin{equation*} 
    \hat{x}_k(t) = \langle e_\alpha, \hat f(t)\rangle \approx (\hat{\mathbf f}(t))^* \mathbf u_\alpha
    \end{equation*}    
    where the vector $\mathbf u_\alpha$ has components $\langle e_\alpha,\psi_\gamma \rangle$, with $\alpha_k = 1$, $\alpha _j = 0$ if $j\neq k$, and $0<|\gamma|\leq d$.
    \end{algorithmic}
\end{algorithm}
This numerical method is reminiscent of the technique developed in \cite{sur16} using the so-called \textit{Koopman observer form}. The main difference is that the Koopman observer form relies on a decomposition onto principal eigenfunctions, while our method considers a decomposition onto $\beta$-unstable eigenfunctions. In the former, (possibly slow) modes associated with non-principal eigenfunctions are not estimated. In the latter, all $\beta$-unstable modes are estimated, which ensures a decay rate greater than $\beta$.
\begin{rem}
    In the last step of Algorithm \ref{alg:obs}, the estimate $\hat{\mathbf x}_k=\langle e_\alpha, \hat{f}(t) \rangle$, with $\vert \alpha \vert = 1$, can be interpreted as the first moment (expectation)  of a distribution. Higher-order moments of this estimated distribution could also be obtained similarly, by computing $\langle e_\alpha, \hat{f}(t) \rangle$ with $\vert \alpha \vert >1$.
\end{rem}

\subsection{Particular class of systems}
We consider a vector field of the form
\begin{equation} \label{dyn_form}
    F_i(\mathbf z) =  -a_i \left( z_i -u_i(\hat{\mathbf{z}}_i)\right) \quad i=1,\ldots, n,
\end{equation}
where $a_i$ are positive real numbers, $u_i\colon \overline{\mathbb D^{n-1}}\to \overline{\mathbb D}$ are holomorphic, and $\hat z_i = (z_1,\cdots, z_{i-1},z_{i+1},\cdots, z_n)$. It follows from Proposition \ref{chelthm1.7} with $G_i(\mathbf z) = - a_i$ that Assumptions 1 and 2 hold if $\text{Re} (\overline{z_i} u_i(\hat{\mathbf z}_i))< 1$ for all $i=1,\ldots,n$.
Clearly, this condition is satisfied if
\begin{equation} \label{hicond}
    \vert u_i(\mathbf w)\vert < 1 \quad \forall \mathbf w \in \mathbb{D}^{n-1}.
\end{equation}
Moreover, if $\partial u_i/\partial z_j(\mathbf 0)=0$ for all $i\neq j$, the eigenvalues of the Jacobian matrix $J_F(\mathbf 0)$ are given by $-a_i$, so that Assumptions 3 and 4 are satisfied provided that the positive real numbers $a_i$ be nonresonant.

\subsection{Numerical experiments}
\begin{figure}            
        \centering
        \includegraphics[width = \linewidth]{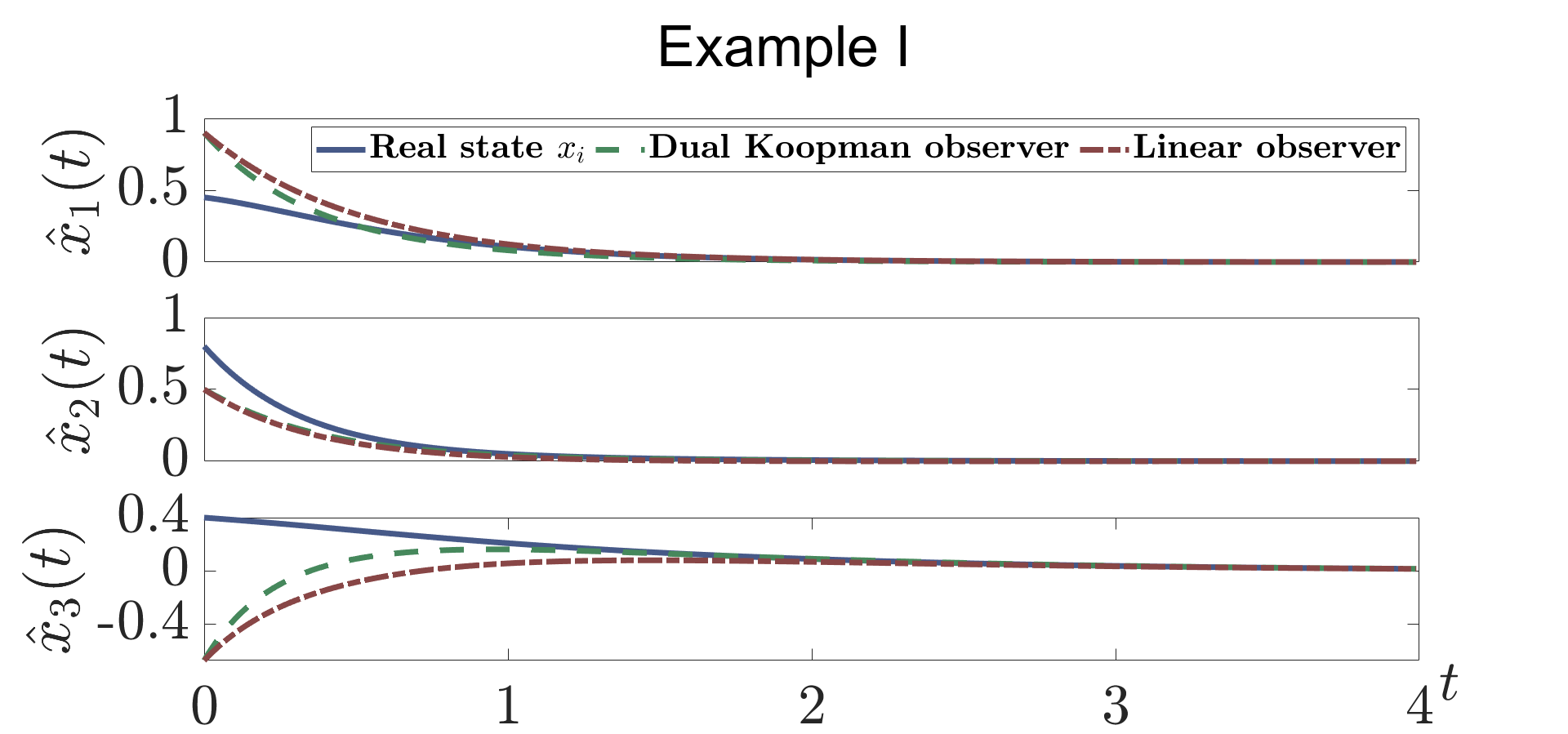}
        \caption{State estimation for the nonlinear system \eqref{exp1}. The initial condition is $\mathbf{x}_0 =(
            0.45, 0.8, 0.4)$ and the initial guess is $\mathbf{\hat{x}}_0 = (0.9, 0.5, -0.67)$.}
        \label{fig:simexp1}
    \end{figure}

\begin{figure}
        \centering
        \includegraphics[width = \linewidth]{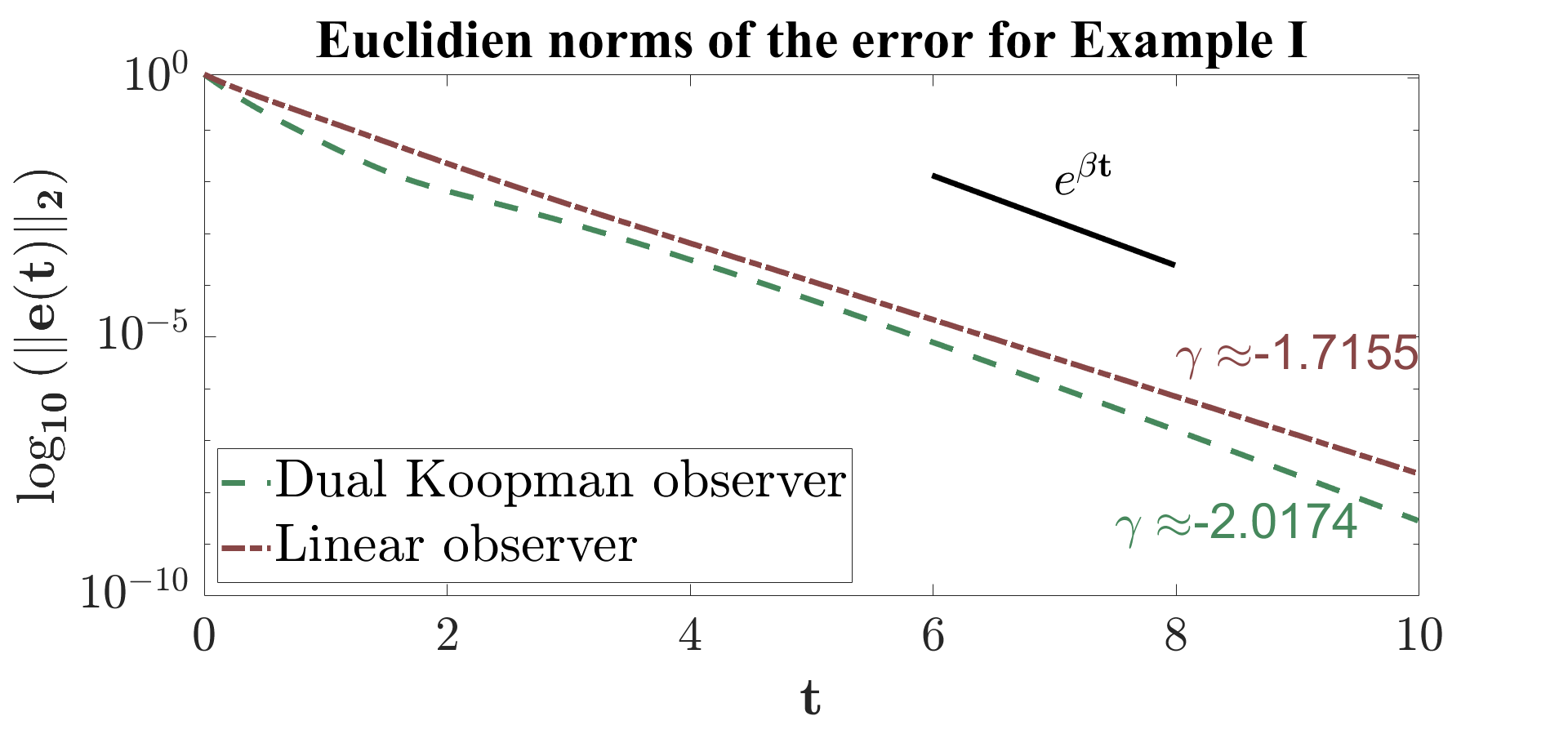}
        \caption{Estimation errors for the nonlinear system \eqref{exp1}.}
        \label{fig:norm}
    \end{figure}

I) We first consider the dynamics
    \begin{equation} \label{exp1}
        \left\{ \begin{array}{l}
            \dot x_1 = - 1.9796~(x_1 - x_2 x_3)\\
            \dot x_2 = -2.95~(x_2 - 0.5x_3^2)\\
            \dot x_3 = -0.853~(x_3 - x_1^2).
        \end{array}
        \right.
    \end{equation}
    It has the form \eqref{dyn_form}, where $u_1(z_2,z_3)=z_2 z_3, u_2(z_1,z_3)=z_3^2/2$ and $u_3(z_1,z_2)=z_1^2$. 
    Those clearly satisfy the condition \eqref{hicond}. Moreover, the parameters have been chosen so that the eigenvalues of the Jacobian matrix $\lambda_1 = -0.853, \lambda_2 = -1.9796$, and $\lambda_3 = -2.95$ are stable and far enough from resonance. The nonlinear output is given by 
    \begin{equation*}
        y(t) = \begin{pmatrix}
            h_1(\mathbf x)\\h_2(\mathbf x) 
        \end{pmatrix} = \begin{pmatrix}
            x_1 + x_2 \\ \cos(x_1) + x_2^2+ x_3.
        \end{pmatrix}.
    \end{equation*}
    and satisfies the condition \eqref{cvg_cond} with $\beta=-2$ ($N_\beta=3$). Recall that the output in the observer will be $y(t)-\begin{pmatrix}
        0 & 1
    \end{pmatrix}^T.$ In Figure \ref{fig:simexp1}, we compare our Koopman observer-based estimate with the one obtained with a Luenberger observer designed for the linearized dynamics. We set the eigenvalues of $(\mathcal A^*_F)^+  + \mathcal L^+ \mathcal C^+$ and the eigenvalues of the closed loop Luenberger observer at $-2, -2.1, -2.2$. The truncation degree is $d=4$. 
    The Koopman observer outperforms the observer based on the linearization during the transient phase, showing a much faster convergence to the true state. In particular, the asymptotic rate of convergence ($\gamma \approx -2.0174$) of the Koopman observer is consistent with the theoretical result (Figure \ref{fig:norm}), and is faster than the rate of convergence of the linearized Luenberger observer ($\gamma \approx -1.7155$), which is approximately equal to $2\lambda_1$. This is explained by the fact that, the linearized observer moves only the modes associated to the principal eigenvalues $\lambda_1$ and $\lambda_2$. Hence, the dominant eigenvalue in closed loop becomes $2\lambda_1$.

\begin{figure}
                   \centering
        \includegraphics[width = \linewidth]{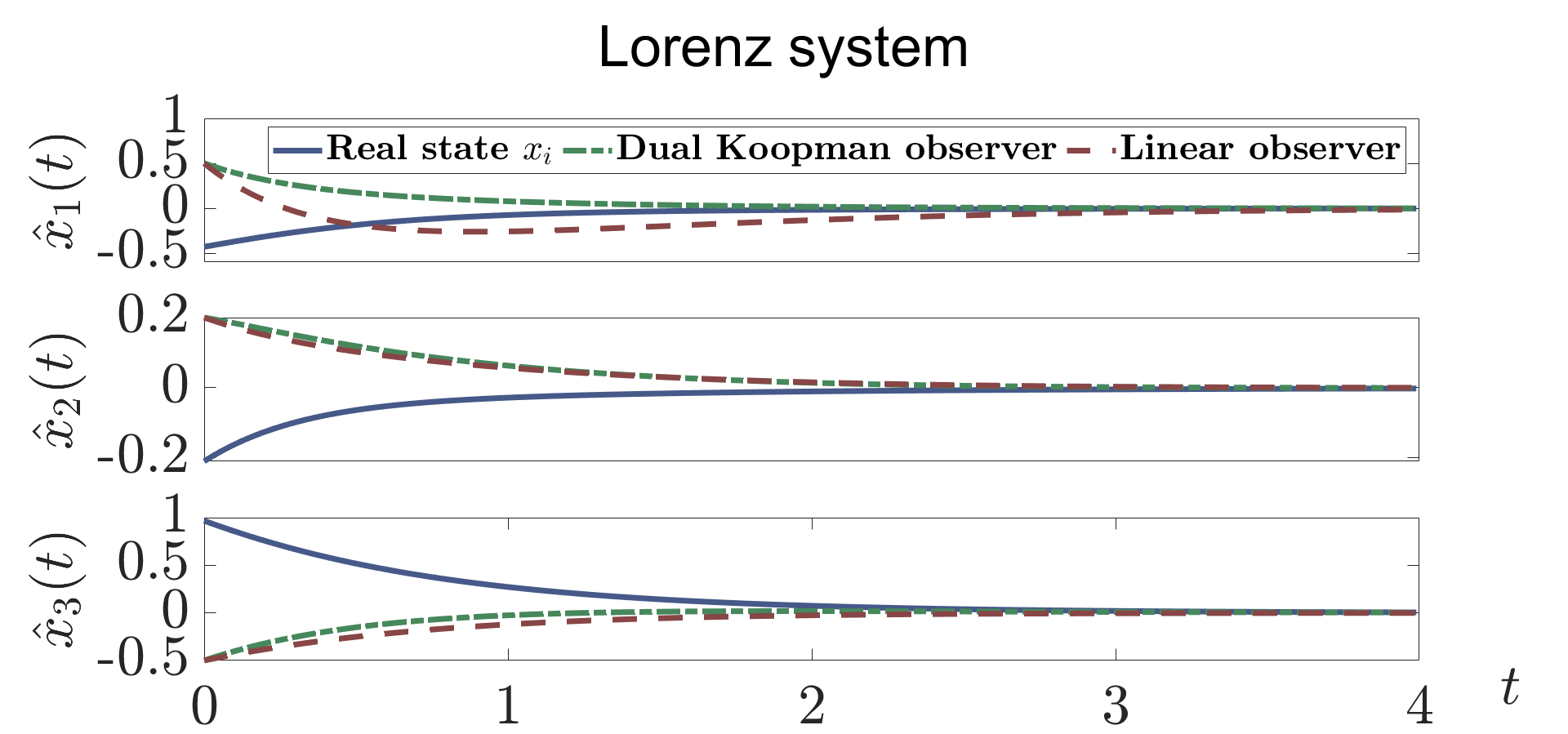}
        \caption{State estimation for the Lorenz system \eqref{ex_lorenz}. The initial condition is $\hat{\mathbf{x}}_0 =(-0.43, -0.21, 0.967)$ and the initial guess is $\hat{\mathbf{x}}_0 = (0.5, 0.2, -0.5)$.} 
         \label{fig:lorenz}
    \end{figure}

\begin{figure}
            
        \centering
        \includegraphics[width = \linewidth]{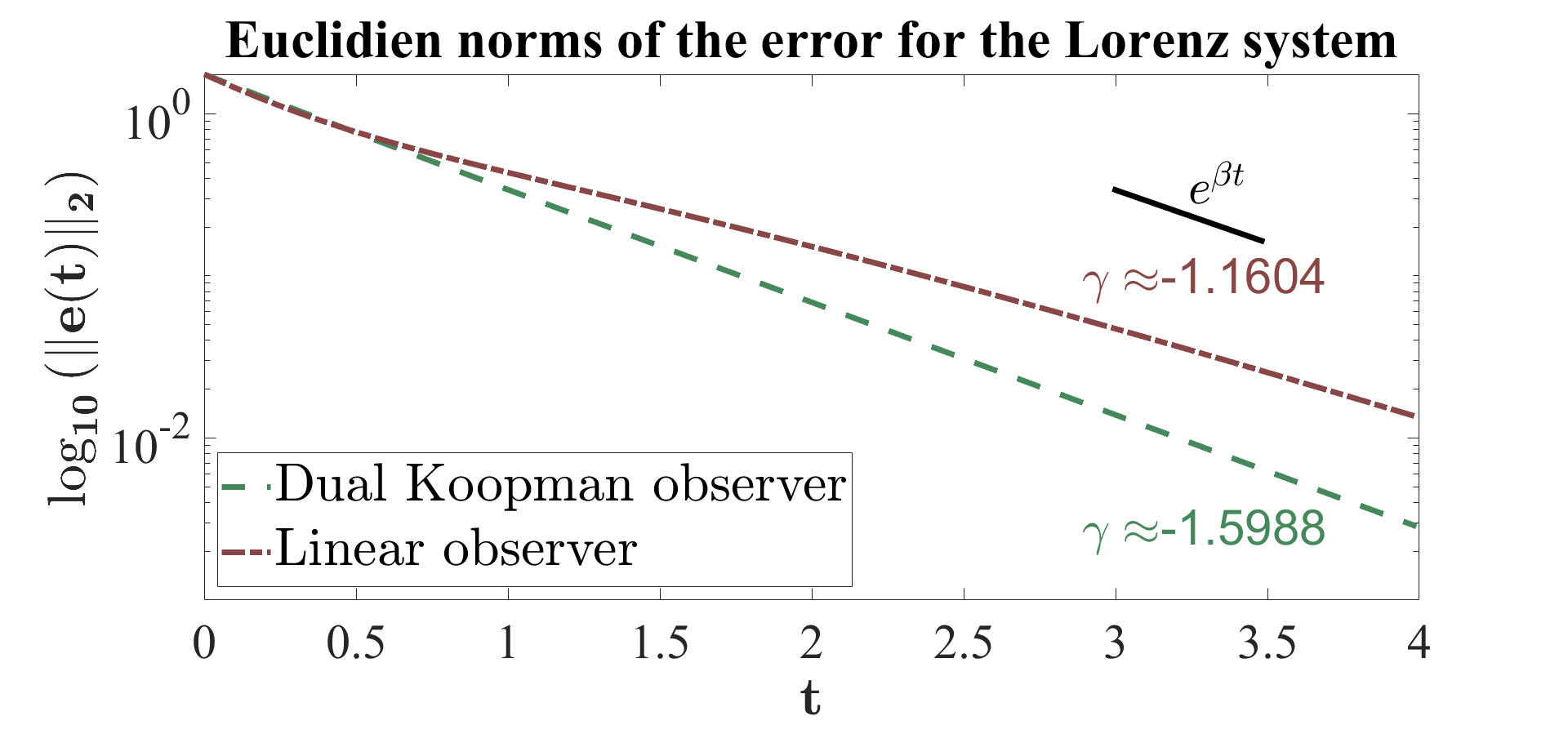}
        \caption{Estimation errors for the Lorenz system \eqref{ex_lorenz}.}
        \label{fig:lorenz_norm}
    \end{figure}

   II) As a prospective example, we consider the Lorenz system
     \begin{equation} \label{ex_lorenz}
        \left\{ \begin{array}{l}
            \dot x_1 = -a x_1 +a x_2\\
            \dot x_2 = rx_1 - x_2 - x_1 x_3\\
            \dot x_3 = -bx_3 +x_1x_2,
        \end{array}
        \right.
    \end{equation}
    with $a= 2.75, b=1.3$ and $r=0.17.$ 
   The eigenvalues of the Jacobian matrix are $\lambda_1 = -0.7645, \lambda_2 = -1.3$, and $\lambda_3 = -2.9855.$ 
    The nonlinear output is given by  
     \begin{equation*}
        y(t) = \begin{pmatrix}
            h_1(\mathbf x)\\h_2(\mathbf x) 
        \end{pmatrix} = \begin{pmatrix}
            \cos(x_1) + x_2 \\ x_1+ x_3.
        \end{pmatrix}.
        \end{equation*}
        and satisfies \eqref{cvg_cond} with $\beta = -1.5$ ($N_\beta=2$). Again, in the observer, the output will be $y(t) - \begin{pmatrix}
            1 & 0
        \end{pmatrix}^T.$ We set the eigenvalues of $(\mathcal A^*_F)^+  + \mathcal L^+ \mathcal C^+$ at $-1.5, -1.6$ and the eigenvalues of the linear Luenberger observer at $-1.5, -1.6, -1.5291$. The truncation degree is $d=6$.
       As shown in Figures \ref{fig:lorenz} and \ref{fig:lorenz_norm}, the Koopman observer-based method outperforms the linear method. It is in particular characterized by a faster rate of convergence $\gamma \approx -1.5988 < \beta$.
       
      Finally, choosing a rate that is too high, i.e. $\beta$ too negative, exposes to large overshoots that are detrimental to quality. This gives raise to a common trade-off between speed of convergence and an adequate input, that could be investigated in further numerical studies.
    
    \section{Conclusion and perspectives}\label{Conclusion}

    This paper sets the foundation for a rigorous theory of the Koopman observer. We have shown that the dual approach is key to state estimation through the Koopman operator framework, and that the RKHS structure plays a crucial role in this context. Moreover, we departed from the strong continuity property and developed the concept of pseudo-weak solution, which allows to leverage the Koopman operator in specific RKHS such as the Hardy space of analytic functions.

  Observability and detectability of the nonlinear system was shown to be equivalent to the so-called pointwise approximate observability and detectability of the dual Koopman system, whereas only sufficient conditions are available in the literature. In the Hardy space on the polydisc, we further derived spectral criteria for observability and detectability, which are based on the principal eigenfunctions of the Koopman operator.

   Our main contribution is to provide a Luenberger-type observer for the dual Koopman system associated with stable hyperbolic equilibrium dynamics. We prove that the estimation error converges weakly and exponentially fast at an arbitrary rate provided that some spectral conditions are satisfied, a property which allows estimating the state in the transient phase, far from the equilibrium.
    
    This present work opens the way to several research perspectives. The class of dynamics considered here is quite limited and could potentially be enlarged to encompass unstable systems and other types of attractors. This could be done by considering other functional spaces associated with weaker regularity properties or other types of linear infinite-dimensional observers. Moreover, while this paper mostly provides a theoretical framework for Koopman operator-based state estimation, the potential of the associated numerical scheme could be further investigated and developed. In particular, finite-dimensional approximation errors could be thoroughly studied, and the numerical method could be extended, for instance, to allow a probabilistic treatment of the initial guess.

\begin{ack}
The authors gratefully acknowledge Christian Zagabe Mugisho (Tu Dortmund) for helpful discussions related to the study of well-posedness. They also thank François-Grégoire Bierwart (UNamur) for his help in the numerical experiments, where they used his toolbox KOSTA (https://github.com/FgBierwart/KOSTA-Toolbox) for specific calculations.
\end{ack}

\section*{Appendix}

\textbf{Proof of Proposition \ref{chelthm1.7}}
Assumption 1 follows from \cite[Theorem 1.7]{che16}. 
Next, observe that, for all $\mathbf z \in \overline{\mathbb D^n}$ and all $t\geqslant 0$, 
    \begin{align}\label{flowderivative} 
        \frac{d}{dt} \vert \varphi^t_j(\mathbf z)\vert^2 &= \frac{d}{dt}\left(\varphi^t_j(\mathbf z) \overline{\varphi^t_j(\mathbf z)}\right) \nonumber\\ 
        &= 2 \,\text{Re}\left( \overline{\varphi^t_j(\mathbf z)} \frac{d}{dt} \varphi^t_j(\mathbf z)\right) \nonumber\\ 
        &= 2 \,\text{Re}\left( \overline{\varphi^t_j(\mathbf z)} F_j(\varphi^t(\mathbf z))\right).
    \end{align} 
    Moreover, for $\mathbf z\in \overline{\mathbb D^n}$ such that $\vert z_j\vert =1$, for $j\in\{1,\ldots,n\}$, we have that 
    \begin{equation*}
        \overline{z_j}F_j(\mathbf z) = (1-\overline{z_j}u_j(\hat{\mathbf z}_j))G_j(\mathbf z).
    \end{equation*}
    Hence, by a uniform continuity argument, it follows from (\ref{holo_flow_ineq}) that there exists $r\in(0,1)$ such that, for all $j\in\{1,\ldots,n\}$  and for all $\mathbf z \in \overline{\mathbb D^n}$ with $r < \vert z_j \vert \leq 1$,
    \begin{equation*} 
        \text{Re } (\overline{z_j} F_j(\mathbf z)) <0.
    \end{equation*}
  Therefore, if there are $t_0 >0$ and $\mathbf z\in \mathbb T^n$ such that $\vert \varphi^{t_0}_j(\mathbf z)\vert \in (r,1)$ for some $j\in\{1,\ldots,n\}$, then, in view of (\ref{flowderivative}),
    the function $t \mapsto \vert \varphi^t_j(\mathbf z)\vert^2$ is strictly decreasing on some time interval $[ t_0,\tau )$. Since the flow is holomorphic, it follows that, for all $t>0$, there exists $R_t < 1$ such that
    \begin{equation*}
      \max_{j\in\{1,\ldots,n\}} \max_{\omega \in [0,2\pi]} \vert \varphi_j^t(i\omega)\vert < R_t.
    \end{equation*}
    Therefore, Assumption 2 holds.
    \hfill\ensuremath{\blacksquare}

\textbf{Proof of Proposition \ref{point_spectrum}}

\emph{Proof of (2).} Recall that 
    \begin{equation*}
        \left\{ \sum_{j=1}^n \alpha_j \lambda_j \colon \alpha \in\mathbb N^n\right\}  \subseteq\sigma\left(A_F\right).
    \end{equation*}
In addition, the principal eigenfunctions of $A_F$ belong to $\mathbb H^2(\mathbb{D}^n)$ (Proposition \ref{principal_ef}), so that the eigenfunctions $\Pi_{j=1}^n \phi_{j}^{\alpha_j}$ associated with the eigenvalues $\sum_{j=1}^n \alpha_j \lambda_j$ also belong to $\mathbb H^2(\mathbb{D}^n)$ for all $\alpha \in \mathbb{N}^n$.
    
Let us show that $\{\psi_\alpha\}_{\alpha\in\mathbb N^n}$ is a basis of $\mathbb H^2(\mathbb{D}^n)$. For \mbox{$d\in\mathbb N$,} we define the matrix $[\mathcal A_F^*]_{d}$ as the upper block triangular submatrix of $\mathcal A_F^*$ containing the rows and columns related to total degrees smaller or equal to $d$.
We denote, by $\mathbf{w}_\alpha^d$ and $\mu_\alpha^d$, $|\alpha|\leq d$ the right eigenvectors and eigenvalues of $[\mathcal A_F^*]_{d}$, respectively.
We verify that, by construction, the infinite-dimensional vector $\mathbf{w}_\alpha = (\mathbf w_\alpha^d \, 0 \, \cdots)^T$ is an eigenvector of $\mathcal A_F^*$ for $|\alpha| \leq d$, and therefore corresponds to the components of $\psi_\alpha$ in the monomial basis. It is also clear that $\psi_\alpha \in \mathbb H^2(\mathbb{D}^n)$.

Next, we show by induction that the eigenfunctions $\psi_\alpha$ are linearly independent. In the case $d=0$ ($\alpha  = 0$), \mbox{$\mathbf{w}_\mathbf{0}=(1 \, 0 \, \cdots)^T$} forms a trivial set of independent vector.
Now, suppose that there is 
$d\in\mathbb N_0$ such that $\{\mathbf w_\alpha\}_{\vert \alpha\vert \leqslant d}$ are linearly independent, and consider
\begin{equation}\label{linindep}
    \sum_{\vert \alpha \vert \leqslant d+1} a_\alpha \mathbf w_\alpha = \mathbf 0
\end{equation}
for some $a_\alpha$. The components of \eqref{linindep} related to the total degree $|\alpha|=d+1$ read $\sum_{\vert \alpha \vert = d+1} a_\alpha \mathbf{\overline{w}}_\alpha^{d+1} = \mathbf 0$ where $\mathbf{\overline{w}}_\alpha^{d+1}$ are the right eigenvectors of the diagonal matrix block associated with the degree $d+1$. These eigenvectors are linearly independent since they are associated with distinct eigenvalues, thanks to the non-resonance assumption (Assumption 4). Hence, $a_\alpha = 0$ if $\vert \alpha \vert = d+1$ and \eqref{linindep} becomes $\sum_{\vert \alpha \vert \leqslant d} a_\alpha \mathbf w_\alpha = 0$. This leads to the conclusion using the induction assumption.

It remains to show that the span of the eigenfunctions $\psi_\alpha$ is dense in $\mathbb H^2(\mathbb{D}^n)$. Pick any $g\in\mathbb H^2(\mathbb{D}^n)$ such that \mbox{$\langle g, \psi_\alpha\rangle = 0$} for all $\alpha\in\mathbb N^n$. 
It follows that $(\mathbf{w}^d_\alpha)^* \mathbf{g}_\alpha^d=0$, where $\mathbf{g}_\alpha^d$ is the vector containing the components of $g$ in the monomial basis up to degree $d$. Since the vectors $\mathbf w_\alpha^d$ are linearly independent, it follows that $\mathbf{g}_\alpha^d=\mathbf 0$ for all $\alpha\in\mathbb N^n$, so that $g=0$. This implies that the span of $\{\psi_\alpha\}_{\alpha\in\mathbb N^n}$ is dense in $\mathbb H^2(\mathbb{D}^n)$, so that $\{\psi_\alpha\}_{\alpha\in\mathbb N^n}$ is a (Schauder) basis of $\mathbb H^2(\mathbb{D}^n)$.

Let us prove that $\{\phi_\alpha\}_{\alpha\in\mathbb N^n}$ and $\{\psi_\beta\}_{\beta\in\mathbb N^n}$ are biorthonormal sequences, up to some rescaling. We just need to show that $\phi_\alpha$ is orthogonal to $\psi_\beta$ for all $\alpha \neq \beta$. We have $\langle \phi_\alpha, \psi_\beta \rangle = (\mathbf w_\beta^d)^* \mathbf v_\alpha^d$, with $d=|\beta|$ and where $\mathbf v^d_\alpha$ and $\mathbf w^d_\beta$ are the vectors of components of $\phi_\alpha$ and $\psi_\beta$, respectively, in the monomial basis up to degree $d$.  Given the block triangular form of $\mathcal{A}_F^*$, the vector $\mathbf v^d_\alpha$ is zero if $|\alpha|>d$ and is the left eigenvector of $[\mathcal{A}_F^*]_d$ otherwise. Recall also that $\mathbf w^d_\beta$ is a right eigenvector of $[\mathcal{A}_F^*]_d$. Hence, $\langle \phi_\alpha, \psi_\alpha \rangle \neq 0$ for all $\alpha$, and $\langle \phi_\alpha, \psi_\beta \rangle = (\mathbf w_\beta^d)^* \mathbf v_\alpha^d=0$ if $\alpha \neq \beta$.

\emph{Proof of (3).}
In order to prove \eqref{expansions}, take any $f\in D(A_F^*)$. Then, 
\begin{equation*}
    A_F^* f = \sum_{\alpha\in\mathbb N^n} \langle A_F^*f, \phi_\alpha\rangle \psi_\alpha
    = \sum_{\alpha\in\mathbb N^n} \overline{\lambda_\alpha} \langle f, \phi_\alpha\rangle \psi_\alpha.
\end{equation*}
In addition, since $K(t)^*$ is a bounded linear operator, we can write, for $f = \sum_{\alpha\in \mathbb N^n} \langle f,\phi_\alpha \rangle \psi_\alpha$,
\begin{equation*}
    K(t)^* f = \sum_{\alpha\in \mathbb N^n} \langle f,\phi_\alpha \rangle K(t)^* \psi_\alpha 
    = \sum_{\alpha\in \mathbb N^n} e^{\overline{\lambda_\alpha} t}\langle f,\phi_\alpha \rangle \psi_\alpha.
\end{equation*}
Similar arguments lead to the representation of $A_F$ and $K(t)$.

\emph{Proof of (1).} It remains to be shown that 
    \begin{equation*}
       \sigma\left(A_F\right) \subseteq S:=  \left\{ \sum_{j=1}^n \alpha_j \lambda_j \colon \alpha \in\mathbb N^n\right\}.
    \end{equation*}

    For this purpose, we adapt the proof of \cite[Theorem 3.2.8]{zwa20}. For all $\lambda \notin S, $
we define the operator 
$$T_{\lambda} f = \sum_{\alpha\in\mathbb N^n} \frac{1}{\lambda - \lambda_\alpha} \langle f,\psi_\alpha\rangle \phi_\alpha .$$
Observe that, for all $f\in\mathbb H^2(\mathbb{D}^n)$, $$\Vert T_\lambda f \Vert \leqslant\displaystyle \sup_{\alpha \in \mathbb N^n} (\vert \lambda - \lambda_\alpha \vert)^{-1} \Vert f\Vert ,$$
where the supremum is finite since the set $S$ is closed. Moreover, using the expansions \eqref{expansions}, we can check that, for all $f\in \mathbb H^2(\mathbb D^n)$ and all $g\in D(A_F)$,
\begin{equation*}
\begin{split}
&(\lambda I - A_F) T_\lambda f = f\\
&T_\lambda (\lambda I -A_F)g = g.    
\end{split} 
\end{equation*}
Hence, the operator $T_\lambda = (\lambda I -A_F)^{-1}$ is linear and bounded on $\mathbb H^2(\mathbb D^n))$ and then, $\lambda$ is in the resolvent set of $A_F$.\hfill\ensuremath{\blacksquare}

\bibliographystyle{apalike}  
\bibliography{source} 

\end{document}